\documentclass[12pt]{amsart} 
\usepackage{amssymb}
\usepackage{latexsym}
\usepackage{amsmath}
\usepackage{amsthm}
\usepackage{amsfonts}
\usepackage{mathrsfs} 
\usepackage{soul}
\usepackage[usenames,dvipsnames]{xcolor}

\usepackage{doi}
\usepackage{tikz-cd} 
\usepackage{lineno} 
\usepackage{datetime} 

\usepackage{ulem}  
\usepackage{cancel} 

\newtheorem{theorem}{Theorem}[section]

\newtheorem{lemma}[theorem]{Lemma}
\newtheorem{corollary}[theorem]{Corollary}
\newtheorem{proposition}[theorem]{Proposition}

\theoremstyle{definition}

\theoremstyle{remark}
\newtheorem{remark}[theorem]{Remark}

\numberwithin{equation}{section}
\numberwithin{theorem}{section}
\allowdisplaybreaks
\usepackage{color}

\newcommand{\R}{\mathbb R}
\newcommand{\N}{\mathbb N}
\newcommand{\C}{\mathbb C}
\newcommand{\T}{\mathbb T}
\newcommand{\Z}{\mathbb Z}

\newcommand{\logl}{L\textnormal{log} L}
\newcommand{\logld}{L(\textnormal{log} L)^2}

\newcommand{\lexp}{L_{\textnormal{exp}}}
\newcommand{\logla}{L(\textnormal{log } L)^{\alpha}}
\newcommand{\loglb}{L(\textnormal{log } L)^{\beta}}
\newcommand{\tlog}{T_{\textnormal{log}}}

\title[The finite   Hilbert transform on  $\logl$]
{The finite   Hilbert transform acting \\ in the Zygmund  space $\logl$}

\author[G.P.  Curbera]{Guillermo P. Curbera}
\address{Facultad de Matem\'aticas \& IMUS,
Universidad de Sevilla, 
Calle Tarfia s/n,  Sevilla 41012, Spain}
\email{curbera@us.es}

\author[S. Okada]{Susumu Okada}
\address{112 Marcorni Crescent, Kambah, ACT 2902, Australia}
\email{susbobby@grapevine.com.au}

\author[W.J. Ricker]{Werner J. Ricker}
\address{Math.--Geogr.\  Fakult\"at, Katholische Universit\"at
Eichst\"att--Ingolstadt, D--85072 Eichst\"att, Germany}
\email{werner.ricker@ku.de}

\date{\today}

\subjclass[2010]{Primary 44A15, 46E30; Secondary  47A53, 47B34.}
\keywords{Finite Hilbert transform, airfoil equation, Zygmund space $\logl$.}

\begin{document}


\begin{abstract} 
The finite Hilbert transform $T$ is a singular integral operator which maps
the Zygmund space $\logl:=\logl(-1,1)$ continuously into $L^1:=L^1(-1,1)$.
By extending the  Parseval and Poincar\'e-Bertrand
formulae to this setting, it is possible to establish an inversion result needed
for solving the airfoil equation  $T(f)=g$ whenever the data function
$g$ lies  in  the range of $T$ within $L^1$ (shown to contain $\logld$). 
Until now this was only known for $g$ 
belonging to the union of all $L^p$ spaces with $p>1$.
It is established  (due to a result of Stein) that $T$ cannot be 
extended to any domain space beyond $\logl$ whilst still taking 
its values in $L^1$, i.e.,  $T\colon \logl\to L^1$ is  optimally defined. 
\end{abstract}

\maketitle


\section{Introduction}
\label{S1}


The finite Hilbert transform $T(f)$ of $f\in L^1:=L^1(-1,1)$ is the  principal value integral
\begin{equation}\label{eq-1-1}
T(f)(t):=\lim_{\varepsilon\to0^+} \frac{1}{\pi}
\left(\int_{-1}^{t-\varepsilon}+\int_{t+\varepsilon}^1\right) \frac{f(x)}{x-t}\,dx ,
\end{equation}
which exists for a.e.\ $t\in(-1,1)$ and is a measurable function.  
It  has important applications to aerodynamics  and elasticity
via the  airfoil  equation 
\begin{equation*} 
\frac{1}{\pi}\, \mathrm{ p.v.}  \int_{-1}^{1}\frac{f(x)}{x-t}\,dx=g(t),
\quad \mathrm{a.e. }\; t\in(-1,1),
\end{equation*}
\cite{cheng-rott}, \cite{muskhelishvili}, \cite{reissner}, \cite{tricomi}, and
to problems arising in image reconstruction; 
see, for example, \cite{katsevich-tovbis}, \cite{sidky-etal}.
We also mention  \cite{duduchava}, \cite{gakhov}, \cite{gohberg-krupnik-1},
\cite{gohberg-krupnik-2},  \cite{mikhlin-prossdorf},
where one-dimensional singular integral operators 
closely related to the finite Hilbert transform are studied in-depth.

In  \cite[Ch.11]{king}, \cite{okada-elliott}, \cite[\S4.3]{tricomi} a detailed investigation
of the finite Hilbert transform was carried out for $T$ acting in the spaces 
$L^p:=L^p(-1,1)$ whenever $1<p<\infty$. 
In \cite{curbera-okada-ricker-ampa},   \cite{curbera-okada-ricker-qm},
\cite{curbera-okada-ricker-mh}, \cite{curbera-okada-ricker-am}
a study of  $T$  
was undertaken when it acts on 
rearrangement invariant (r.i.)  spaces 
$X$ on $(-1,1)$ having non-trivial Boyd indices $0<\underline{\alpha}_X\le \overline{\alpha}_X<1$.
This larger class of spaces provides a natural setting for a further study
of $T$ beyond the $L^p$-spaces, as is exemplified by  two facts. First,  that 
$T\colon X\to X$ is injective if and only if $L^{2,\infty}\not\subseteq X$ 
and second, that $T\colon X\to X$ has a 
non-dense range if and only if $X \subseteq L^{2,1}$  (for $X$ separable). 
Here $L^{2,1}$ and $L^{2,\infty}$ are 
the usual Lorentz spaces, which are r.i.\    spaces but not  $L^p$-spaces.

In all investigations so far $T$ is always considered as a linear operator acting
from a r.i.\  space  into \textit{itself}. The point of
departure of this paper is to consider $T$ acting in the classical
Zygmund space  $\logl:=\logl(-1,1)$. This is a r.i.\  space
on $(-1,1)$ close to $L^1$ in the sense that
$L^p\subseteq\logl$ for all $1<p<\infty$, which implies that 
$X\subseteq \logl$ for all r.i.\   spaces 
$X$  with $0<\underline{\alpha}_X\le \overline{\alpha}_X<1$. 
But,  the Boyd indices of $\logl$ are trivial 
(both equal 1) and so $T$ \textit{cannot} map $\logl$ into itself;
see Section 2. However, it turns out that $T$ \textit{does} map
$\logl$ continuously into the strictly larger space $L^1$
(cf.\  Theorem \ref{t-2-1}). 
The aim of this paper is to make an in-depth study of the operator $T\colon\logl\to L^1$.

Section \ref{S2} is mainly devoted to presenting some preliminary facts 
concerning Banach function spaces, r.i.\ spaces 
and their Boyd indices, the Zygmund space $\logl$, the finite Hilbert transform $T$ and
its action in $\logl$.

Fundamental tools for studying  $T$ are the Parseval formula
$$
\int_{-1}^{1}fT(g)=-\int_{-1}^{1}gT(f),
$$
and the Poincar\'e-Bertrand  (or   the Hardy-Poincar\'e-Bertrand, \cite[p.57]{king}) formula
$$
T(gT(f)+fT(g))=(T(f))(T(g))-fg,\quad \mathrm{a.e.},
$$
for suitable pairs of functions $f, g$.
Tricomi, \cite[\S4.3 (2), (4)]{tricomi}, and 
Love, \cite[Corollary]{love}, showed that these formulae
hold whenever $f\in L^p$ and $g\in L^{p'}$, with $\frac1p+\frac1{p'}=1$. 
The validity  of these formulae was extended to the setting of r.i.\ spaces
$X$, showing that they
hold for all $f$ in $X$ and $g$ in the associate space of $X$, provided that 
$0<\underline{\alpha}_X\le \overline{\alpha}_X<1$,
\cite[Proposition 3.1]{curbera-okada-ricker-ampa}.
In Section \ref{S3} suitable extensions of these formulae,
beyond the class  $\bigcup_{1<p< \infty} L^p
=\bigcup_{0<\underline{\alpha}_X\le \overline{\alpha}_X<1} X$ and
needed in this paper for the space  $\logl$, 
are established  in Theorems \ref{t-3-1} and \ref{t-3-2}.
These theorems provide the means to determine the kernel of the operator
$T\colon\logl\to  L^1$. The arcsine distribution 
$x\mapsto1/\sqrt{1-x^2}$, for $x\in(-1,1)$,
plays an important role in the study of $T$. The reason is that 
\begin{equation*}
T\Big(\frac{1}{\sqrt{1-x^2}}\Big)(t)=
\text{p.v.} \frac{1}{\pi}\int_{-1}^{1}\frac{1}{\sqrt{1-x^2} \,(x-t)}\,dx=0,\quad t\in(-1,1),
\end{equation*}
and, moreover, as  proved by Tricomi,   
if $T(f)(t)=0$ for a.e.\ $t\in(-1,1)$ with 
$f\in\bigcup_{1<p<\infty} L^p$, then necessarily $f$ is a constant
multiple of $1/\sqrt{1-x^2}$, \cite[\S4.3 (14)]{tricomi}. 
For a recent extension of 
this result we refer to \cite{coifman-steinerberger};  also \cite{huang} is related to this
extension. 
Tricomi's result is extended by showing that the kernel of 
$T\colon\logl\to  L^1$  is the one-dimensional
space spanned by 
$1/\sqrt{1-x^2}$ (cf.\ Theorem \ref{t-3-4}).

As alluded to above, the arcsine distribution 
pervades the entire theory of the finite Hilbert transform.
For instance, it occurs  in the definition of an important auxiliary
operator $\widehat T$, given by
\begin{equation*}
\widehat{T}(g)(x):=\frac{-1}{\sqrt{1-x^2}}\, 
T\big(\sqrt{1-t^2}\,g(t)\big)(x),\quad \mathrm{a.e. }\; x\in (-1,1),
\end{equation*}
which plays a central role in considerations of the inversion 
of $T$ in the $L^p$-spaces, for $1<p<2$, where it is
known that $\widehat T$ is bounded from $L^p$ into itself.  
In Section 4, via extrapolation theorems of Yano, it is shown that both
$T\colon L(\textnormal{log } L)^{1+\beta}\to L(\textnormal{log } L)^{\beta}$
and $\widehat{T}\colon L(\textnormal{log } L)^{1+\beta}\to L(\textnormal{log } L)^{\beta}$
are bounded operators for all $\beta\ge0$ (cf.\ 
Propositions \ref{p-4-4} and \ref{p-4-5}, respectively).
Various features connecting  $T$ and $\widehat T$ in the space $\logl$
are presented in Theorem \ref{t-4-10}. In Corollary \ref{cor-4-10} a description
of the range space $T(\logl)$ is given. Both of these results are crucial for establishing
an inversion formula needed for solving the airfoil equation.
Tricomi studied the solution of this equation via the Poincar\'e-Bertrand formula 
and proved an inversion formula   within the class of $L^p$-spaces for $1<p<\infty$
with $p\not=2$, \cite[Ch.11]{king}, \cite{okada-elliott}, \cite[\S4.3 (10)--(12)]{tricomi}.
The solution of the airfoil equation was extended to  
$T\colon X\to X$, with $X$ 
in the larger class of all r.i.\  spaces  whose Boyd indices satisfy
$0<\underline{\alpha}_X\le \overline{\alpha}_X<1/2$ or 
$1/2<\underline{\alpha}_X\le \overline{\alpha}_X<1$,
\cite[Corollary 3.5]{curbera-okada-ricker-ampa}.
In the former case, $T\colon X\to X$ is injective 
and has a closed range with co-dimension 1 and 
an explicit formula for the range is possible,
\cite[Theorem  3.3]{curbera-okada-ricker-ampa}, while in the latter
case $T\colon X\to X$  has a one-dimensional  kernel and  is surjective,
\cite[Theorem 3.2]{curbera-okada-ricker-ampa}.
The solution of the airfoil equation $T(f)=g$,
where $g$ is a given function in the range of $T\colon\logl\to L^1$
and $f\in\logl$ is to  be found, is presented in Theorem \ref{t-4-13}. 
Note  that the range space $T(\logl)$ contains $\logld$
(cf. Proposition \ref{p-4-11}(ii)), which properly contains 
the union of all r.i.\ spaces $X$ with non-trivial Boyd indices.
Subtleties concerning various properties of the range space
$T(\logl)\subseteq L^1$ are exposed in Proposition \ref{p-4-11}.

Let $\T$ denote the circle group, $F\colon L^1(\T)\to c_0(\Z)$ be the Fourier transform and 
$2\le p'\le \infty$. The following question was raised by R.\ E.\ Edwards. 
\textit{What can be said about the family $[F,\ell^{p'}(\mathbb{Z})]$
consisting of all functions $f\in L^1(\mathbb{T})$ having the property that 
$F(f\chi_A)\in \ell^{p'}(\mathbb{Z})$ for all Borel sets $A\subseteq\mathbb{T}$?},  \cite[p.206]{edwards}.
The requirement that $F(f\chi_A)\in \ell^{p'}(\mathbb{Z})$ for all Borel sets 
$A\subseteq\mathbb{T}$ distinguishes $[F,\ell^{p'}(\mathbb{Z})]$ from the linear space
$F^{-1}(\ell^{p'}(\mathbb{Z}))$. Indeed, it has the consequence that 
$[F,\ell^{p'}(\mathbb{Z})]$ is an \textit{ideal} of functions in $L^1(\T)$.
The Hausdorff-Young inequality ensures that $L^p(\T)\subseteq [F,\ell^{p'}(\mathbb{Z})]$.
The precise identification of $[F,\ell^{p'}(\mathbb{Z})]$
is given in \cite{mockenhaupt-ricker}, where it is shown that
$L^p(\T)\subsetneqq [F,\ell^{p'}(\mathbb{Z})]$, that is, 
the Fourier transform operator $F\colon L^p(\T)\to \ell^{p'}(\Z)$
is \textit{not} optimally defined but has an extension to the genuinely larger
Banach function space $[F,\ell^{p'}(\mathbb{Z})]$ and this space
is optimal; no further extension of $F$ (with values still in $\ell^{p'}(\Z)$) is
possible. Given   a r.i.\ space $X$  over $(-1,1)$, the same problem as above  arises.
\textit{Identify, if possible, the family $[T,X]$ of all functions
$f\in L^1$ such that $T(f\chi_A)\in X$ for all Borel sets $A\subseteq(-1,1)$}.
Whenever $X$ has non-trivial Boyd indices this was achieved in
\cite{curbera-okada-ricker-ampa}, \cite{curbera-okada-ricker-mh},
where it was shown that $[T,X]=X$, that is, 
$T\colon X\to X$ is already optimally defined. With $[T,L^1]$ 
denoting the space of all functions $f\in L^1$ such that
$T(f\chi_A)\in L^1$ for all Borel sets $A\subseteq(-1,1)$, it is
clear that $X\subseteq [T,L^1]$ for every r.i.\ space with non-trivial indices. 
In Section \ref{S5}   (due to a result of Stein) we identify $[T,L^1]$ 
precisely as the classical Zygmund space $\logl$; see Theorem \ref{t-5-6}. 
Consequently, the continuous linear operator 
$T\colon\logl \to L^1$  is optimally defined.  
Noticeable properties of $T$ needed to establish
Theorem \ref{t-5-6}, some giving equivalent conditions for
a function to belong to  $\logl$ in terms
of $T$, are derived in Proposition \ref{p-5-1}.

In order to keep the presentation of the main results as transparent
as possible we have  placed some technical results   in the final 
Section \ref{S6} (an Appendix).


The first author acknowledges the support  of 
PID2021-124332NB-C21 FEDER/Ministerio de Ciencia e Innovaci\'on and 
FQM-262 (Spain).


\section{Preliminaries}
\label{S2}


The setting of this paper is the measure space  $(-1,1)$ 
equipped with its Borel $\sigma$-algebra $\mathcal{B}$ and  Lebesgue measure $|\cdot|$
(restricted to $\mathcal{B}$). Denote by 
$L^0:=L^0(-1,1)$ the space (of equivalence classes) of all $\mathbb{C}$-valued
measurable functions, endowed with the topology of convergence in measure. 
The space $L^p(-1,1)$ is denoted simply by $L^p$, for $1\le p\le\infty$.

A \textit{Banach function space} (B.f.s.) $X$ on  $(-1,1)$ is a
Banach space  $X\subseteq L^0$ satisfying
the ideal property, that is, $g\in X$ and $\|g\|_X\le\|f\|_X$
whenever $f\in X$ and $|g|\le|f|$ a.e.  
The associate space $X'$  of $X$ consists  of all
functions $g$ satisfying $\int_{-1}^1|fg|<\infty$, for every
$f\in X$, equipped with the norm
$\|g\|_{X'}:=\sup\{|\int_{-1}^1fg|:\|f\|_X\le1\}$. 
The space $X'$ is a closed subspace of the dual Banach space $X^*$ of $X$. 
If $f\in X$ and $g\in X'$, then $fg\in L^1$ and
$\|fg\|_{L^1}\le \|f\|_X \|g\|_{X'}$, i.e., H\"older's inequality is available. 
The second  associate space $X''$ 
of $X$ is defined as $X''=(X')'$. The norm in $X$ is called absolutely continuous if,
for every $f\in X$, we have $\|f\chi_A\|_X\to0$ whenever $|A|\to0$.
The space $X$ satisfies the Fatou property  if, 
whenever  $(f_n)_{n=1}^\infty\subseteq X$ satisfies
$0\le f_n\le f_{n+1}\uparrow f$ a.e.\ with $\sup_n\|f_n\|_X<\infty$,
then $f\in X$ and $\|f_n\|_X\to\|f\|_X$.   
In this paper \textit{all} B.f.s.' $X$, as in 
\cite{bennett-sharpley},    
are assumed to satisfy the Fatou property. In this case $X''=X$ and hence,
$f\in X$ if and only if $\int_{-1}^1|fg|<\infty$, for every $g\in X'$.
Moreover, $X'$ is a norm-fundamental subspace of $X^*$, that is, 
$\|f\|_X=\sup_{\|g\|_{X'}\le1} |\int_{-1}^1fg|$ for $f\in X$. 
If $X$ is separable, then $X'=X^*$.

A \textit{rearrangement invariant} (r.i.) space $X$ on $(-1,1)$ is a B.f.s.\  
such that if $g^*\le f^*$ with $f\in X$,  
then $g\in X$ and $\|g\|_X\le\|f\|_X$.
Here $f^*\colon[0,2]\to[0,\infty]$ is 
the decreasing rearrangement of $f$, that is, the
right continuous inverse of its distribution function
$\lambda\mapsto|\{t\in (-1,1):\,|f(t)|>\lambda\}|$ for $\lambda\ge0$.
The associate space $X'$ of a r.i.\ space $X$ is again a r.i.\ space.
Every r.i.\ space on $(-1,1)$ satisfies 
$L^\infty\subseteq X\subseteq L^1$.

The family of r.i.\ spaces includes many classical spaces 
appearing in analysis, such as  the Lorentz $L^{p,q}$ spaces, 
\cite[Definition IV.4.1]{bennett-sharpley}, Orlicz $L^\varphi$ spaces 
\cite[\S4.8]{bennett-sharpley}, Marcinkiewicz $M_\varphi$ spaces, 
\cite[Definition II.5.7]{bennett-sharpley}, Lorentz $\Lambda_\varphi$ spaces,
\cite[Definition II.5.12]{bennett-sharpley},
and the Zygmund $L^p(\text{log L})^\alpha$ spaces, 
\cite[Definition IV.6.11]{bennett-sharpley}. In particular, 
$L^p=L^{p,p}$, for $1\le p\le \infty$.  
The space weak-$L^1$, denoted by $L^{1,\infty}:=L^{1,\infty}(-1,1)$, will play
an important role; it is  a quasi-Banach space
and satisfies $L^1\subseteq L^{1,\infty} \subseteq L^0$, with both
inclusions continuous.

The Zygmund space $\logl:=\logl(-1,1)$ consists of all measurable functions $f$
on $(-1,1)$ for which either one of the following two equivalent conditions hold:
$$
\int_{-1}^1|f(x)|\log^+|f(x)|\,dx<\infty,\quad \int_0^2 f^*(t)\log\Big(\frac{2e}{t}\Big)\,dt<\infty,
$$
see \cite[Definition IV.6.1 and Lemma IV.6.2]{bennett-sharpley}. 
The space $\logl$ is  r.i.\  with absolutely continuous norm
(cf. \cite[p.247]{bennett-sharpley}) given by
$$
\|f\|_{\logl}:=\int_0^2 f^*(t)\log\Big(\frac{2e}{t}\Big)\,dt,\quad f\in\logl.
$$
It follows that $L^p\subseteq \logl$ for all $1<p<\infty$.  The
associate space of $\logl$ is the space $\lexp$ consisting 
of all measurable functions $f$ on $(-1,1)$ having exponential integrability;
see \cite[Definition IV.6.1]{bennett-sharpley}. 
The separability of  $\logl$ implies that $(\logl)^*=(\logl)'=\lexp$.

We will also require the family of spaces $\logla$ for $\alpha>1$,
consisting of all measurable functions $f$
on $(-1,1)$ for which either one of the following two equivalent conditions hold:
$$
\int_{-1}^1|f(x)|(\log(2+|f(x)|)^\alpha\,dx<\infty,\quad \int_0^2 f^*(t)\log^\alpha
\Big(\frac{2e}{t}\Big)\,dt<\infty.
$$
The space $\logla$ is  r.i.\  with absolutely continuous norm given by
$$
\|f\|_{\logla}:=\int_0^2 f^*(t)\log^\alpha\Big(\frac{2e}{t}\Big)\,dt,\quad f\in\logla;
$$
see \cite[Definition IV.6.11 and Lemma IV.6.12]{bennett-sharpley}.  
The following inclusions hold: 
$$
L^p\subseteq \loglb\subseteq \logla\subseteq \logl,
\quad 1<p,\,\, 1<\alpha<\beta.
$$

Given a r.i.\ space $X$ on $(-1,1)$,  the Luxemburg representation theorem 
ensures that there exists a r.i.\  space  $\widetilde X$ on $(0,2)$
such that $\|f\|_X=\|f^*\|_{\widetilde X}$ for $f\in X$, \cite[Theorem II.4.10]{bennett-sharpley}.
The dilation operator $E_t$ for $t>0$ is defined, for 
each $f\in \widetilde X$, by $E_t(f)(s):=f(st)$ for $0\le s\le \min\{2,1/t\}$ and zero 
for  $\min\{2,1/t\}< s\le 2$. The operator $E_t\colon \widetilde X\to \widetilde X$  is bounded 
with $\|E_{1/t}\|_{\widetilde X\to \widetilde X}\le \max\{t,1\}$,
\cite[Theorem III.5.11]{bennett-sharpley}. The \textit{lower} and \textit{upper 
Boyd indices} of  $X$ are defined, respectively, by
\begin{equation*}
\underline{\alpha}_X\,:=\,\sup_{0<t<1}\frac{\log \|E_{1/t}\|_{\widetilde X\to \widetilde X}}{\log t}
\;\;\mbox{and}\;\;
\overline{\alpha}_X\,:=\,\inf_{1<t<\infty}\frac{\log \|E_{1/t}\|_{\widetilde X\to \widetilde X}}{\log t} ;
\end{equation*}
see  \cite[Definition III.5.12]{bennett-sharpley}. 
They satisfy $0\le\underline{\alpha}_X\le \overline{\alpha}_X\le1$.
Note that $\underline{\alpha}_{L^p}= \overline{\alpha}_{L^p}=1/p$
for all $1\le p<\infty$, and $\underline{\alpha}_{\logl}= \overline{\alpha}_{\logl}=1$,
\cite[Theorem IV.6.5]{bennett-sharpley}.


For  a r.i.\ space $X$,  the boundedness of the finite Hilbert  transform
$T\colon X\to X$ (which is indicated by simply writing $T_X$)
is equivalent to $X$ having non-trivial  Boyd indices, that is, 
$0<\underline{\alpha}_X\le \overline{\alpha}_X<1$, 
see  \cite[pp.170--171]{krein-petunin-semenov}. 
Since $\underline{\alpha}_{X'}=1-\overline{\alpha}_X$
and $\overline{\alpha}_{X'}=1-\underline{\alpha}_X$, the condition 
$0<\underline{\alpha}_X\le \overline{\alpha}_X<1$ implies that 
$0<\underline{\alpha}_{X'} \le \overline{\alpha}_{X'}<1$. Hence,
 $T_{X'}\colon X'\to X'$ is also bounded. 
It follows  from the Parseval formula for $T_X$, 
(cf.\  \cite[Proposition 3.1(b)]{curbera-okada-ricker-ampa}) 
that the restriction of the dual operator 
$T_X^*\colon X^*\to X^*$ of $T_X$ to the closed subspace $X'$ of $X^*$ 
is precisely $-T_{X'}\colon X'\to X'$.

The operator $T$ is not continuous on $L^1$. However, 
due to a result of Kolmogorov $T\colon L^1\to L^{1,\infty}$ is continuous, 
\cite[Theorem III.4.9(b)]{bennett-sharpley}. 
The  following result is central to this paper.


\begin{theorem}\label{t-2-1}
The finite Hilbert transform $T\colon\logl\to L^1$ is a bounded operator.
\end{theorem}

\begin{proof}
The  Calderon operator 
is defined for every measurable function $f$ on $(0,\infty)$ by 
$$
S(f)(t):=\int_0^\infty\min\Big(1,\frac{s}{t}\Big)f(s)\,ds=
\frac{1}{t}\int_0^tf(s)\,ds+\int_t^\infty f(s)\frac{ds}{s},\quad t>0;
$$
see (4.22) of \cite[\S III.4, p.133]{bennett-sharpley}. Clearly 
$0\le S(g)\le S(f)$ whenever $0\le g\le f$. We will apply $S$ to the
decreasing rearrangement $f^*$ of a function $f$.
For the Hilbert transform $H$ on $\R$ it is known that there exists
a  constant $c>0$ satisfying 
\begin{equation}\label{eq-2-1}
(H(f))^*(t)\le c\cdot S(f^*)(t),\quad t>0, 
\end{equation}
for all functions $f$ satisfying the condition $S(f^*)(1)<\infty$,
\cite[Theorem III.4.8]{bennett-sharpley}.

Given $f\in L^1$, note that $S(f^*)(1)<\infty$ 
is satisfied because $f^*\colon(0,2)\to \R$ and
\begin{align*}
S(f^*)(1)&=
\int_0^1f^*(s)\,ds+\int_1^2 f^*(s)\frac{ds}{s}
\le \int_0^2f^*(s)\,ds=\|f\|_{L^1}<\infty.
\end{align*}
Observe, for $0<t<2$, that
\begin{equation}\label{eq-2-3}
S(f^*)(t)=
\frac{1}{t}\int_0^tf^*(s)\,ds+\int_t^2 f^*(s)\frac{ds}{s} .
\end{equation}
Moreover, since $f\in L^1$, it follows that $T(f)=-\chi_{(-1,1)}H(f\chi_{(-1,1)})$.
Then, for each $0<t<2$, it follows from \eqref{eq-2-1} and \eqref{eq-2-3} that
\begin{align*}
(T(f))^*(t)&=\big(-\chi_{(-1,1)}H(f\chi_{(-1,1)})\big)^*(t)
\le \big(H(f\chi_{(-1,1)})\big)^*(t) 
\\ & \le c\cdot S(f\chi_{(-1,1)})^*(t) \le  c\cdot S(f^*)(t)
= c\cdot\frac{1}{t}\int_0^tf^*(s)\,ds+c\int_t^2 f^*(s)\frac{ds}{s}.
\end{align*}
Integrating  this inequality and applying Fubini's theorem  yields
\begin{align*}
\|T(f)\|_{L^1}&=\|(T(f))^*\|_{L^1(0,2)}=\int_0^2(T(f))^*(t)\,dt
\\ & \le c\int_0^2\frac{1}{t}\int_0^tf^*(s)\,ds\,dt+c\int_0^2\int_t^2 f^*(s)\frac{ds}{s}\,dt 
\\ & = c\int_0^2f^*(s)\log\Big(\frac{2}{s}\Big)\,ds+c\int_0^2f^*(s)\,ds 
\\ & \le c\int_0^2f^*(s)\log\Big(\frac{2e}{s}\Big)\,ds+c\int_0^2f^*(s)\,ds 
\\ & = c\|f\|_{\logl}+c\|f\|_{L^1}
\le 2c\|f\|_{\logl}.
\end{align*}
So, whenever $f\in \logl$ we can conclude that
$\|T(f)\|_{L^1}\le 2c\,\|f\|_{\logl}$.  
Accordingly, $T\colon\logl\to L^1$ is a bounded operator (with $\|T\|\le 2c$).
\end{proof}


For all of the above (and further) facts on r.i.\  spaces see \cite{bennett-sharpley}, 
\cite{krein-petunin-semenov}, \cite{lindenstrauss-tzafriri}, for example.

The definition of $T$ in \eqref{eq-1-1} coincides with that used by  
Tricomi, \cite[\S4.3]{tricomi}, whereas King uses $-T$, \cite[Ch.11]{king},
and both J\"orgens, \cite[\S13.6]{jorgens}, Widom, \cite{widom}, use $-iT$.
When convenient, the operator $T\colon\logl\to L^1$ will
also be denoted by $\tlog$, where the subscript $\text{log}$ indicates briefly the domain 
space $\logl$.

 
\section{Parseval and Poincar\'e-Bertrand formulae}
\label{S3}


In this section we present  integrability conditions for a pair of 
functions   $f, g$ which ensure the validity of 
the Parseval and the Poincar\'e-Bertrand formulae
beyond the class  $\bigcup_{1<p<\infty} L^p
=\bigcup_{0<\underline{\alpha}_X\le \overline{\alpha}_X<1} X$;
see \cite[Proposition 2.b.3]{lindenstrauss-tzafriri}.
These are then applied to determine the kernel of $\tlog\colon \logl\to L^1$.


\begin{theorem}\label{t-3-1}
Let  the functions $f\in L^1$ and $g\in \logl$ satisfy   $fT(g\chi_A)\in L^1$, 
for every  set $A\in\mathcal B$.
Then  $gT(f)\in L^1$ and the following Parseval formula is valid:
\begin{equation}\label{eq-3-1}
\int_{-1}^{1}fT(g)=-\int_{-1}^{1}gT(f).
\end{equation}
\end{theorem}


\begin{proof}
\textit{Step 1.} Suppose first that $f \in L^\infty$ and $g\in  \logl$, 
in which case the assumption of the theorem is automatically satisfied.  That is,
$fT (g\chi_A) \in L^1$ for every $A\in \mathcal B$ (because $T (g\chi_A) \in L^1$; 
see Theorem \ref{t-2-1}).  Fix $A \in \mathcal B$.  
The aim is to show that $(g\chi_A)T(f) \in L^1$ and
\begin{equation}\label{eq-3-2}
\int_{-1}^{1}fT(g\chi_A)=-\int_{-1}^{1}(g\chi_A)T(f).
\end{equation}
Let $A_n:=|g|^{-1}([0,n])$ and $g_n:=g\chi_{A_n}$ for $n\in\N$.
Then $(g_n\chi_A)T(f)\to (g\chi_A)T( f)$ pointwise.
Given $B \in \mathcal B$, both of the functions $f$ and 
$g_n\chi_{A\cap B}$ belong to $L^\infty$ and hence, also to $L^2$, for every $n\in\N$.  
Applying the Parseval formula for $L^2$ yields
\begin{equation}\label{eq-3-3}
\int_{-1}^1fT(g_n\chi_{A\cap B})=-\int_{-1}^1(g_n\chi_{A\cap B})T(f) 
= - \int_B (g_n\chi_A) T(f),\quad n\in\N.
\end{equation}
Since $g_n\chi_{A\cap B} \to g\chi_{A\cap B}$ in $\logl$, it follows from 
Theorem \ref{t-2-1} that $T(g_n\chi_{A\cap B} )\to T(g\chi_{A\cap B})$ 
in $L^1$.  This and \eqref{eq-3-3} imply that
\begin{equation}\label{eq-3-4}
\int_{-1}^1fT(g\chi_{A\cap B})=\lim_n \int_{-1}^1fT(g_n\chi_{A\cap B})
= - \lim_n\int_B (g_n\chi_A)T(f).
\end{equation}
So, the sequence
$\big(\int_B (g_n\chi_A)T(f)\big)_{n=1}^\infty$  is convergent in $\C$ 
for every  set $B\in\mathcal B$.  Thus, the pointwise limit  function 
$(g\chi_A)T(f) = \lim_n (g_n\chi_A)T(f)$ is integrable and
\begin{equation}\label{eq-3-5}
\lim_n (g_n\chi_A)T(f)= (g\chi_A)T(f)\quad \text{in }L^1.
\end{equation}
This and \eqref{eq-3-4}, with $B:= (-1,1)$, ensure that \eqref{eq-3-2} is valid as
$$
\int_{-1}^1fT(g\chi_A)=-\lim_n\int_{-1}^1(g_n\chi_A)T(f) = -\int_{-1}^1 (g\chi_A)T(f).
$$

\medskip

\textit{Step 2.}  Now let   $f\in L^1$ and   $g\in \logl$   satisfy the assumption of the theorem.
Set $B_n:=|f|^{-1}([0,n])$ and $f_n:=f\chi_{B_n}$ for $n\in\N$, in which case
$f_n \to f$ in $L^1$.  
Apply Kolmogorov's theorem to conclude that $T(f_n) \to T(f)$ in 
measure.  So, we can suppose that $T(f_n) \to T(f)$ 
pointwise a.e., by passing to a subsequence, if necessary.

Fix $A\in \mathcal B$.  Apply Step 1 to $f_n\in L^\infty$ and $g\in \logl$ 
to conclude that $(g\chi_A)T(f_n) \in L^1$ and
\begin{equation}\label{eq-3-6}
\int_{-1}^{1}f_nT(g\chi_A)=-\int_{-1}^{1}(g\chi_A)T(f_n),\quad n\in\N.
\end{equation}
Since $|f_nT(g\chi_A)|\le |fT(g\chi_A)|\in L^1$ for each $n\in\N$,  
dominated convergence  ensures that
\begin{equation}\label{eq-3-7}
\lim_n f_nT(g\chi_A )=  fT(g\chi_A ) \quad \text{in } L^1
\end{equation} 
and that
\begin{equation*}\lim_n \int_{-1}^{1}f_nT(g\chi_A) = -\int_{-1}^{1}fT(g\chi_A).
\end{equation*} 
Hence, it follows from \eqref{eq-3-6} that
\begin{equation}\label{eq-3-8}
\lim_n\int_{A}gT(f_n)=-\lim_n\int_{-1}^{1}f_nT(g\chi_A)=-\int_{-1}^{1}fT(g\chi_A).
\end{equation}
That is,   for every  set $A\in\mathcal B$, the sequence
$\big(\int_AgT(f_n)\big)_{n=1}^\infty$  is convergent in $\C$.
This fact, together with $gT(f_n)\to gT(f)$ a.e.,  imply  that $gT(f)\in L^1$ and
\begin{equation}\label{eq-3-9}
\lim_n gT(f_n)= gT(f) \quad \text{in } L^1.
\end{equation}
From \eqref{eq-3-8},  with $A: = (-1,1)$, and  \eqref{eq-3-9} we can derive
 the Parseval formula \eqref{eq-3-1}.
\end{proof}


\begin{theorem}\label{t-3-2}
Let  the functions $f\in L^1$ and $g\in \logl$ satisfy   
$fT(g\chi_A)\in L^1$, for every  set $A\in\mathcal B$.
Then   the following Poincar\'e-Bertrand  formula 
is valid in $L^0$:
\begin{equation}\label{eq-3-10}
T(gT(f)+fT(g))=(T(f))(T(g))-fg.
\end{equation}
\end{theorem}


\begin{proof}
Theorem \ref{t-3-1} ensures that the functions $f$ and $g$ satisfy both 
$gT(f) \in L^1$ and the Parseval formula \eqref{eq-3-1}.  Let  $B_n:=|f|^{-1}([0,n])$ 
and $f_n:=f\chi_{B_n}$ for $n\in\N$, as in  Step 2 of  the proof of Theorem  \ref{t-3-1}.  Then 
the following conditions are satisfied:
\begin{itemize}
\item[(a)] $T\big(gT(f_n)\big) \to T\big(gT(f)\big)$ in $L^0$,
\item[(b)]  $T\big(f_nT(g)\big) \to T\big(fT(g)\big)$ in $L^0$,
\item[(c)] $T(f_n) T(g)\to  T(f)T(g)$ in $L^0$, and
\item[(d)]  $f_ng\to fg$ pointwise and hence, in $L^0$.
\end{itemize}
Indeed, recall from \eqref{eq-3-9}  that $\lim_n gT(f_n)= gT(f)$ in $L^1$ 
and  from \eqref{eq-3-7}, with $A: =(-1,1)$, that  $\lim_n f_nT(g)=  fT(g)$   in $L^1$.   
So, Kolmogorov's theorem  ensures (a) and (b), respectively.  
Moreover, the convergence $T(f_n)\to T(f) $ in measure 
has already  been  verified in   Step 2 of the 
proof of Theorem \ref{t-3-1}.  Hence, (c) holds.  Condition (d) is clearly satisfied.

Fix $n\in \N$.  The claim is, as an identity in $L^0$, that
\begin{equation}\label{eq-3-11}
T\big(gT(f_n)+f_nT(g)\big)=\big(T(f_n)\big)\big(T(g)\big)-f_ng .
\end{equation}
To verify this  we follow the argument used in 
Step 1 for the proof of Theorem \ref{t-3-1}.  Let 
$A_k: = |g|^{-1}([0,k])$ and $g_k:= g\chi_{A_k}$ for $k\in\N$. 
Then \eqref{eq-3-5}, with $f_n$ in place of $f$ and $A: = (-1,1)$, 
implies that $\lim_k g_kT(f_n) = gT(f_n)$ in $L^1$. Then Kolmogorov's theorem yields
\begin{equation}\label{eq-3-12}
\lim_kT\big( g_kT(f_n)\big) =T\big( gT(f_n)\big) \quad \text{in }L^0.
\end{equation}
Since $g_k\to g$ in $\logl$, it follows from Theorem \ref{t-2-1} 
that $\lim_kT(g_k)= T(g)$ in $L^1$ and hence, that $\lim_k f_nT(g_k) = f_nT(g)$
in $L^1$ because $f_n \in L^\infty$.  Kolmogorov's theorem again implies that
\begin{equation}\label{eq-3-13}
\lim_kT\big( f_nT(g_k)\big) =T\big( f_nT(g)\big) \quad \text{in }L^0.
\end{equation}
Now, the functions $f_n$ and $g_k$ with $k\in\N$ all belong to $L^\infty$ 
and hence, also to $L^2$.  By the  Poincar\'e-Bertrand  formula for $L^2$ we obtain
the identity
\begin{equation}\label{eq-3-14}
T\big(g_kT(f_n)+f_nT(g_k)\big)=\big(T(f_n)\big)\big(T(g_k)\big)-f_ng_k 
\end{equation}
in $L^1$ for every  $k\in \N$.
Therefore,  \eqref{eq-3-11}  follows from  \eqref{eq-3-12}, 
\eqref{eq-3-13} and \eqref{eq-3-14} because
\begin{align*}
T\big(gT(f_n)+f_nT(g)\big)
&=  \lim_k T\big(g_kT(f_n)+ f_nT(g_k)\big)
\\&
=\lim_k\Big(\big (T(f_n)\big)\big(T(g_k)\big)-f_ng_k\Big)
= \big(T(f_n)\big)\big( T(g) \big)- f_ng
\end{align*}
with both limits existing in $L^0$.

Finally we can obtain the Poincar\'e-Bertrand  formula
\eqref{eq-3-10} from \eqref{eq-3-11} by applying (a)--(d) as follows:
\begin{align*}
T\big(gT(f)+fT(g)\big) &=\lim_nT\big(gT(f_n)+f_nT(g)\big)
\\&=\lim_n\Big(\big(T(f_n)\big)\big( T(g) \big)- f_ng\Big)
=\big(T(f)\big)\big(T(g)\big)-fg.
\end{align*}
\end{proof}


\begin{corollary}\label{c-3-3}
The following Parseval  formula 
\begin{equation}\label{eq-3-15}
\int_{-1}^1fT(g)=-\int_{-1}^1gT(f),\quad f\in L^\infty,\; g\in\logl, 
\end{equation}
holds, as does the Poincar\'e-Bertrand formula (in $L^0$) 
\begin{equation}\label{eq-3-16}
T(gT(f)+fT(g))=(T(f))(T(g))-fg,\quad f\in L^\infty,\; g\in\logl.
\end{equation}
\end{corollary}

\begin{proof} 
The Parseval formula \eqref{eq-3-15} was verified in Step 1 of the proof
of Theorem \ref{t-3-1}. Moreover, Theorem \ref{t-3-2} implies that 
 \eqref{eq-3-16} is valid.
\end{proof}


The following theorem extends to $\logl$ Tricomi's identification of  the kernel of $T$ on $L^p$; 
see also   \cite[Theorems 3.2(a) and 3.3(a)]{curbera-okada-ricker-ampa} 
for other r.i.\ spaces $X$. 
Note, via \cite[Theorem IV.6.5]{bennett-sharpley}, that 
$$
\bigcup_{1<p<\infty} L^p
=\bigcup_{0<\underline{\alpha}_X\le \overline{\alpha}_X<1} X
\subsetneqq \logl.
$$


\begin{theorem}\label{t-3-4}
Let $f\in \logl$. Then  $T(f)=0$ in $L^1$
if and only if $f(x)=c/\sqrt{1-x^2}$, for some 
constant $c\in\C$.
\end{theorem}

\begin{proof}
That $T(f)(x)=0$ for a.e.\  $x\in(-1,1)$, 
when $f(x)=c/\sqrt{1-x^2}$, is a straight-forward verification; see \cite[\S4.3 (7)]{tricomi}.

Let $f\in \logl$ satisfy $T(f)(x)=0$  for a.e.\  $x\in(-1,1)$.  Set
$g(x)=\sqrt{1-x^2}\in L^\infty$. Then $gT(f)=0$. Moreover, 
$T(g)(x)=-x$ (see \cite[(11.57)]{king}).
Applying \eqref{eq-3-16}, with the roles of $f$ and $g$ interchanged, yields
$$
T(tf(t))(x)=f(x)\sqrt{1-x^2}.
$$
But, the function $T(tf(t))$ is  constant because,  for a.e.\  $x\in(-1,1)$, we have
\begin{align*}
T(tf(t))(x)&=\frac1\pi \text{p.v.} \int_{-1}^{1}\frac{tf(t)}{t-x}\,dt=
\frac1\pi \int_{-1}^{1}f(t)\,dt
+\frac{x}{\pi}  \text{p.v.} \int_{-1}^{1}\frac{f(t)}{t-x}\,dt
\\ & 
=\frac1\pi \int_{-1}^{1}f(t)\,dt+xT(f)(x)
=\frac1\pi \int_{-1}^{1}f(t)\,dt.
\end{align*}
It follows  that $\sqrt{1-x^2} f(x)$ is  a constant function. 
\end{proof}


\begin{remark}\label{r-3-5} 
Corresponding to important properties of the arcsine distribution
 $1/\sqrt{1-x^2}$ in relation to the operator $T$, is the role played by the 
Marcinkiewicz space $L^{2,\infty}$, 
also known as  weak-$L^2$, \cite[Definition IV.4.1]{bennett-sharpley},
which consists of those functions  $f\in L^0$  satisfying, for some constant $M>0$,
the inequality
\begin{equation*}
f^*(t)\le \frac{M}{t^{1/2}},\quad 0<t\le2.
\end{equation*}
Indeed, $L^{2,\infty}$ is the \textit{smallest} r.i.\ space  
containing the function $1/\sqrt{1-x^2}$
(whose decreasing rearrangement is the function
$t\mapsto 2/t^{1/2}$ for $t\in(0,2)$),  that is, 
$1/\sqrt{1-x^2}\in X$ if and only if $L^{2,\infty}\subseteq X$. 
Hence, as mentioned before, 
$T\colon X\to X$ is injective if and only if $L^{2,\infty}\not\subseteq X$.
\end{remark}


\section{Inversion of the finite Hilbert transform on $\logl$}
\label{S4}


The operator $T$ maps $\logl$ continuously into $L^1$ (cf.\ Theorem \ref{t-2-1}).
We will require the extension of  this result for $T$ acting in  
the larger class of spaces $\logla$ for $\alpha\ge0$.
This is obtained via extrapolation and will be needed for the inversion formula. 
For the case of $T_p\colon L^p\to L^p$ with $1<p<2$
and  the case of $T_X\colon X\to X$ with $1/2<\underline{\alpha}_X\le \overline{\alpha}_X<1$,
the inversion formula  is given in terms of an operator $\widehat T$ (see \eqref{eq-4-1})
which satisfies $\widehat T\colon L^p\to L^p$ boundedly for $1<p<2$. We also need to extend
the boundedness of $\widehat T$ to the class of spaces $\logla$ for $\alpha\ge0$.
As for $T$, this is again obtained via extrapolation.


The following is a classical result on extrapolation; 
see  \cite[Theorem XII.(4.41)]{zygmund}, \cite[Theorem 2.1]{edmunds-krbec}.

\begin{theorem}[Yano]\label{t-4-1}
Let $1<p_0<\infty$ and  $S$ be a linear  operator that 
maps $L^p$ boundedly into $L^p$
for all $1<p<p_0$ and  such that there exist  constants $C>0$ and 
$1<\alpha\le p_0$ satisfying 
$$
\|S\|_{L^p\to L^p}\le \frac{C}{p-1},\quad p\in(1,\alpha).
$$
Then $S$ can be extended  to $\logl$ with $S\colon \logl\to L^1$
a bounded operator.
\end{theorem}

A related result is the following one, \cite[Theorem 5.1]{edmunds-krbec}.

\begin{theorem}\label{t-4-2}
Let $p_0, S, C$ and $\alpha$ satisfy
the conditions in Yano's theorem.
Suppose, for some $\gamma\ge0$, that 
$S\colon L(\textnormal{log } L)^{\gamma}\to L^1$ 
is bounded. Then 
$$
S\colon L(\textnormal{log } L)^{\gamma+\beta}\to L(\textnormal{log } L)^{\beta}   
$$
boundedly, for all $\beta>0$.
\end{theorem}


In order to apply Yano's theorem to  $T$
we will require the following inequality.

\begin{lemma}\label{l-4-3}
Let $1<p<2$. The operator norm of $T\colon L^p\to L^p$ satisfies
\begin{equation*}
\|T\|_{L^p\to L^p}  \le\frac{3}{p-1}.
\end{equation*}
\end{lemma}

\begin{proof}  
For each $1<p<\infty$, Pichorides  proved in \cite{pichorides}  that
the Hilbert transform $H\colon L^p(\R)\to L^p(\R)$
satisfies
$$
\|H\|_{L^p(\R)\to L^p(\R)}=  \max\big\{\tan(\pi/(2p)),\cot(\pi/(2p))\big\},
\quad 1<p<\infty.
$$
We know that $T\colon L^p\to L^p$ boundedly. 
A result of McLean and Elliott, \cite[Theorem 3.4]{mclean-elliott} yields that 
$\|T\|_{L^p\to L^p}=\|H\|_{L^p(\R)\to L^p(\R)}$. So,
\begin{equation*}
\|T\|_{L^p\to L^p} = \max\big\{\tan(\pi/(2p)),\cot(\pi/(2p))\big\},
\quad 1<p<\infty.
\end{equation*}
Consequently, for $1<p<2$, we have $\|T\|_{L^p \to L^p} = \tan(\pi/(2p))$.

Observe  that $\sin\theta \ge \theta/2$ whenever 
$ 0\le \theta \le \pi/3$.  Since 
$\frac{\pi}{2} -\frac{\pi}{2p} =\frac{\pi}{2} 
\big(1 -\frac{1}{p}\big)<\frac{\pi}{2} \big(1 -\frac{1}{2}\big) < \frac{\pi}{3} $, 
it follows that
$$
\cos \Big(\frac{\pi}{2p}\Big) = \sin\Big(\frac{\pi}{2} -\frac{\pi}{2p}\Big) \ge
\frac{1}{2}\Big(\frac{\pi}{2} -\frac{\pi}{2p}\Big) =
\frac{\pi}{4}\Big(1-\frac{1}{p}\Big)=\frac{\pi(p-1)}{4p} . 
$$
This implies  that
$$
\|T\|_{L^p \to L^p} = \tan \Big(\frac{\pi}{2p}\Big) 
= \frac{\sin\big(\frac{\pi}{2p}\big)}{\cos\big(\frac{\pi}{2p}\big)} 
\le \frac{4p}{\pi(p-1)} < \frac{8}{\pi(p-1)} < \frac{3}{p-1}.
$$
\end{proof}


Observe that Theorem  \ref{t-4-1} and Lemma \ref{l-4-3} provide an alternate
(perhaps more abstract) proof of Theorem \ref{t-2-1}. The following 
more general result  follows from  Theorem \ref{t-4-2}, with 
$\alpha=p_0=2$ and $S=T$, $C=3$, $\gamma=1$, together with Lemma \ref{l-4-3}.

\begin{proposition}\label{p-4-4}
The finite Hilbert transform $T$ satisfies
\begin{equation*}
T\colon L(\textnormal{log } L)^{1+\beta}\to L(\textnormal{log } L)^{\beta}
\end{equation*}
boundedly, for every  $\beta\ge0$.
\end{proposition}


Recall, for  $g\in L^1$, that one can define pointwise a measurable function 
$\widehat{T}(g)\in L^0$ by
\begin{equation}\label{eq-4-1}
\widehat{T}(g)(x):=\frac{-1}{\sqrt{1-x^2}}\, 
T\big(\sqrt{1-t^2}\,g(t)\big)(x),\quad \mathrm{a.e. }\; x\in (-1,1).
\end{equation}
It is known, for $1<p<2$,  that  $\widehat{T}\colon L^p\to L^p$ boundedly,
\cite[Theorem I.4.2]{gohberg-krupnik-1}.
In order to apply extrapolation, as was  done for $T$ in Proposition \ref{p-4-4},
we need an explicit upper bound   on the operator norms
$\|\widehat{T}\|_{L^p\to L^p}$ for $p$ near 1, which is not available in \cite{gohberg-krupnik-1}.
To achieve this we will require some auxiliary facts. 
The desired result is as follows.


\begin{proposition}\label{p-4-5}
For each $\beta\ge0$ the operator
\begin{equation*}
\widehat{T}\colon L(\textnormal{log } L)^{1+\beta}\to L(\textnormal{log } L)^{\beta}
\end{equation*}
boundedly. In particular, $\widehat{T}\colon \logl\to L^1$ is bounded.
\end{proposition}


Given numbers $\beta,\gamma \in (0,1)$ satisfying $1 < (\beta + \gamma)$ define
\begin{equation}\label{eq-4-2}
c(\beta,\gamma):=
\max\Big\{\frac{1}{1-\beta},\frac{1}{1-\gamma},\frac{1}{\beta+\gamma-1}
\Big\}.
\end{equation}

The proof of the following technical result is given in the Appendix.

\begin{lemma}\label{l-4-6}
Let $\beta,\gamma \in (0,1)$ satisfy $1 < (\beta + \gamma)$. Then
\begin{equation}\label{eq-4-3}
\int_{-1}^\infty\frac{d\xi}{|\xi|^{\beta}(\xi+1)^{\gamma}}
\le 6\; c(\beta,\gamma)
\end{equation}
and
\begin{equation}\label{eq-4-4}
 \int_{-1}^1\frac{dt}{|t-x|^{\beta}|1-t^2|^{\gamma}}
\le \frac{24\, c(\beta,\gamma)}{(1-x^2)|^{\beta+\gamma-1}},\quad x\in(-1,1).
\end{equation}
\end{lemma}


A consequence of the previous inequalities is the following result.

\begin{lemma}\label{l-4-7}
Let $1<p<2$. For each $\delta$  such that
$$
\frac{(p-1)}{2}<\delta p<\min\Big\{\frac12, (p-1)\Big\}
$$
the norm of the linear operator $(T+\widehat{T})\colon L^p\to L^p$
satisfies
$$
\|T+\widehat{T}\|_{L^p\to L^p}\le 
\frac{24\sqrt2}{\pi}\Big(c\Big(\frac12,\frac12+\delta p\Big)\Big)^{1/p} 
\Big(c\Big(\frac12,\delta p'\Big)\Big)^{1/p'}.
$$ 
\end{lemma}

\begin{proof}
Given $g \in L^p$ and $x \in (-1,1)$ it follows, 
via $|a^{1/2} -b^{1/2}|\le |a-b|^{1/2}$, that

\begin{align*} 
|(T +\widehat{T})(g)(x)|
&=
\frac{1/\pi}{(1-x^2)^{1/2}} \left|  \mathrm{p.v.}\int_{-1}^1 \frac{(1-x^2)^{1/2} g(t)}{t-x} dt
-  \mathrm{p.v.}\int_{-1}^1 \frac{(1-t^2)^{1/2} g(t)}{t-x} dt\right|
\\&\le 
\frac{1/\pi}{(1-x^2)^{1/2}}
\int_{-1}^1 \frac{\big|(1-x^2)^{1/2}- (1-t^2)^{1/2}\big|\cdot |g(t)|}{|t-x|}\,dt
\\ &\le 
\frac{1/\pi}{(1-x^2)^{1/2}}
\int_{-1}^1 \frac{\big|(1-x^2)- (1-t^2)\big|^{1/2}|g(t)|}{|t-x|}\,dt 
\\ &=
\frac{1/\pi}{(1-x^2)^{1/2}}\int_{-1}^1 \frac{|t-x|^{1/2}\,|t+x|^{1/2} |g(t)|}{|t-x|}\,dt
\\ &\le
\frac{\sqrt2/\pi}{(1-x^2)^{1/2}}\int_{-1}^1 \frac{|g(t)|}{|t-x|^{1/2}}dt
\qquad (\text{via}  \quad |t+x|\le 2)
\\ &=
\frac{\sqrt2/\pi}{(1-x^2)^{1/2}}\int_{-1}^1 \frac{(1-t^2)^\delta|g(t)|}{|t-x|^{1/2}(1-t^2)^{\delta}}dt
\\ &\le
\frac{\sqrt2/\pi}{(1-x^2)^{1/2}}\left(\int_{-1}^1 \frac{(1-t^2)^{\delta p}|g(t)|^p}{|t-x|^{1/2}}dt\right)^{1/p}
\left(\int_{-1}^1 \frac{dt}{|t-x|^{1/2}(1-t^2)^{\delta p'}}\right)^{1/p'}.
\end{align*}
The last inequality is via H\"older's inequality.
Now apply \eqref{eq-4-4} with $\beta: =1/2$ and  $\gamma :=\delta p'$ to see 
(as $\beta+\gamma-1=\delta p'-\frac12$)
that
$$
\int_{-1}^1 \frac{1}{|t-x|^{1/2}(1-t^2)^{\delta p'}}dt 
\le \frac{24\, c\big(\frac{1}{2}, \delta p'\big)}{(1-x^2)^{ \delta p' -(1/2)}}.
$$
Hence,
$$
|(T +\widehat{T})(g)(x)|\le 
\frac{ \sqrt2\Big(24\, c\big(\frac{1}{2}, \delta p'\big)\Big)^{1/p'}}{ \pi(1-x^2)^{\delta +(1/(2p))}}
\left(\int_{-1}^1 \frac{(1-t^2)^{\delta p}|g(t)|^p}{|t-x|^{1/2}}dt\right)^{1/p}.
$$
This   yields, via Fubini's theorem, that
\begin{align*}
\|(T +\widehat{T})(g)\|_{L^p}^p &\le \Big(\frac{\sqrt2}{\pi}\Big)^p 
\Big(24\, c\big(\frac{1}{2}, \delta p'\big)\Big)^{p/p'} 
\int_{-1}^1\frac{1}{(1-x^2)^{ \delta p +(1/2)}}\left(\int_{-1}^1 
\frac{(1-t^2)^{\delta p}|g(t)|^p}{|t-x|^{1/2}}dt\right)dx
\\ &=
\Big(\frac{\sqrt2}{\pi}\Big)^p\Big(24\, c\big(\frac{1}{2}, \delta p'\big)\Big)^{p/p'}
\int_{-1}^1 \frac{|g(t)|^p}{(1-t^2)^{-\delta p}}\left(\int_{-1}^1 
\frac{dx}{|t-x|^{1/2}(1-x^2)^{ \delta p +(1/2)}}  \right)dt
\\&\le 
\Big(\frac{\sqrt2}{\pi}\Big)^p \Big(24\, c\big(\frac{1}{2}, \delta p'\big)\Big)^{p/p'}
\Big( 24\,c\big( \frac{1}{2}, \delta p + \frac{1}{2}\big)\Big) \int_{-1}^1 |g(t)|^pdt
\end{align*}
because  \eqref{eq-4-4}, with $\beta := 1/2$ and 
$\gamma:= \delta p +(1/2)\in(0,1)$, ensures that
$$
\int_{-1}^1 \frac{1}{|t-x|^{1/2}(1-x^2)^{\delta p+ (1/2)}} dx 
\le \frac{24\,c\big(\frac{1}{2},  \delta p +\frac{1}{2}\big)}{(1-t^2)^{1/2 +( \delta p + (1/2)) -1}} =
\frac{24\,c\big(\frac{1}{2},  \delta p +\frac{1}{2}\big)}{(1-t^2)^{ \delta p }}.
$$
Since  $g \in L^p$ is arbitrary, it follows that
$$
\|T +\widehat{T}\|_{L^p\to L^p} \le \frac{\sqrt2}{\pi}
\cdot 24 \left(c\Big(\frac{1}{2}, \delta p + \frac{1}{2}\Big)\right)^{1/p}
\left(c\Big (\frac{1}{2} ,\delta p' \Big)\right)^{1/p'}.
$$  
\end{proof}


\begin{lemma}\label{l-4-8}
Let $1<p<3/2$ and set $\delta:=\frac23\cdot \frac{(p-1)}{p}$.   
The operator norm of $T+\widehat{T}$ satisfies
$$
\|T+\widehat{T}\|_{L^p\to L^p}\le 
\frac{24\sqrt2}{\pi}\Big(\frac{3}{p-1}\Big)^{1/p} 
6^{1/p'}\le \frac{72\sqrt2}{p-1}.
$$ 
\end{lemma}

\begin{proof}
Note that $\min\{\frac12, p-1\}=p-1$ and that 
$\frac12(p-1)<\delta p<(p-1)$, as $\delta p=\frac23(p-1)$.
Moreover, $\delta$ satisfies 
$\frac12(p-1)<\delta p<\min\big\{\frac12, (p-1)\big\}$.
It then follows from Lemma \ref{l-4-7}  and Lemma \ref{l-6-1}   that
$$
\|T +\widehat{T}\|_{L^p \to L^p}
\le \frac{24\sqrt2}{\pi}
\left (\frac{3}{p-1}\right)^{1/p} \cdot 6^{1/p'}
\le \frac{72\sqrt2}{p-1}.
$$
\end{proof}


It is now  possible to establish Proposition \ref{p-4-5}.
\medskip

\noindent
\textit{Proof of Proposition \ref{p-4-5}.} 
Let $1<p<\frac32$. It follows from Lemma \ref{l-4-3} and Lemma \ref{l-4-8} that
$$
\|\widehat{T}\|_{L^p\to L^p}\le 
\|T+\widehat{T}\|_{L^p\to L^p} + \|T\|_{L^p\to L^p}
\le \frac{72\sqrt2}{p-1}+ \frac{3}{p-1}.
$$
This proves Proposition \ref{p-4-5} by applying Theorem \ref{t-4-2} with $S: =T$. 
\qed


\begin{remark}\label{r-4-9}
The following diagram illustrates the action of  $T$ and $\widehat{T}$ 
in various spaces.
\begin{center}
\begin{tikzcd}
L(\textnormal{log } L)^{3} \arrow[r, hook, ] \arrow[dr,  "T"] 
& L(\textnormal{log } L)^{2} \arrow[r, hook, ] \arrow[dr,  "T"] 
& L(\textnormal{log } L) \arrow[dr,  "T"]  
& 
\\
& L(\textnormal{log } L)^{2} \arrow[r, hook, ] 
& L(\textnormal{log } L) \arrow[r, hook, ] 
&  L^1
\\
L(\textnormal{log } L)^{3} \arrow[r, hook, ] \arrow[ur,  "\widehat{T}"]
& L(\textnormal{log } L)^{2} \arrow[r, hook, ] \arrow[ur,  "\widehat{T}"]
& L(\textnormal{log } L) \arrow[ur,  "\widehat{T}"] 
&  
\end{tikzcd}
\end{center}
\end{remark}


We are now in a suitable situation for  studying the action of  $T$ and 
$\widehat{T}$ on $\logl$.

\begin{theorem}\label{t-4-10}  
The following assertions are valid.
\begin{itemize}
\item[(i)]   The operator $\tlog\colon\logl\to L^1$ is not injective.
\item[(ii)] Let  $g\in\logl$.   Then $\widehat T(g)\in L^1$ and $T(\widehat T(g))=g$. Moreover,
\begin{equation*}
\int_{-1}^1\widehat{T}(g)(x)\,dx=0.
\end{equation*}
\item[(iii)] The operator $P\colon\logl\to\logl $  given by
\begin{equation*}
P (f)(x):=\left(\frac1\pi\int_{-1}^1 f(t)\,dt\right)
\frac{1}{\sqrt{1-x^2}},  \quad |x|<1, \quad f\in {\logl} ,
\end{equation*}
is a bounded projection satisfying 
$$
\|P\|_{\logl\to\logl}\le \frac1\pi\Big\|\frac{1}{\sqrt{1-x^2}}\Big\|_{\logl}.
$$
Furthermore, $TP\colon\logl\to L^1$ is the zero operator.
\item[(iv)] For each $f\in\logl$ it is the case that 
\begin{equation}\label{eq-4-5}
f-P(f)=\widehat{T}(T(f)).
\end{equation}
Moreover, $\widehat{T}T\colon\logl\to\logl$ satisfies
\begin{equation}\label{eq-4-6}
\|\widehat{T}T\|_{\logl\to\logl}
\le 1+\frac1\pi\Big\|\frac{1}{\sqrt{1-x^2}}\Big\|_{\logl}.
\end{equation}
\end{itemize}
\end{theorem}


\begin{proof} 
(i) Theorem \ref{t-3-4} shows that
$\tlog^{-1}(\{0\}) = \mathrm{span} \{1/\sqrt{1-x^2}\}\subseteq \logl$.

(ii) Let $w(x):=\sqrt{1-x^2}$ for $x\in(-1,1)$. Since $|g\chi_A|\le|g|$, 
Proposition \ref{p-4-5} shows that
$\widehat{T}(g\chi_A) \in L^1$ for every $A \in \mathcal B$ and so,
via \eqref{eq-4-1}, also $(1/w)T(wg\chi_A)= -\widehat{T}(g\chi_A) \in L^1$.
Then Theorem \ref{t-3-1}, with $f:=1/w$ and $(gw)$ in place of $g$, implies  that 
$(gw)T(1/w) \in L^1$ (as already known since $T(1/w) =0$)  and,
via \eqref{eq-3-1},  that also $\int_{-1}^1(1/w)T(gw) =-\int_{-1}^1 gwT(1/w) = 0$.  
This identity yields $\int_{-1}^1 \widehat{T}(g)(x)\,dx = 0$.
Finally, apply \eqref{eq-3-10} in Theorem \ref{t-3-2}, still with $f: =1/w$ 
and $(gw)$ in place of $g$, to obtain 
$$
T(\widehat{T}(g))= -T\Big((1/w)T(gw)\Big) = T\Big(gwT(1/w)\Big)
-T(1/w)T(gw) +(1/w)(gw) =g.
$$

(iii) The boundedness of $P$ follows from
$$
\|P(f)\|_{\logl }\le \frac1\pi \left\|\frac{1}{\sqrt{1-x^2}}\right\|_{\logl } 
\|f\|_{L^1}
\le \frac1\pi \left\|\frac{1}{\sqrt{1-x^2}}\right\|_{\logl } 
\|f\|_{\logl },
$$
which in turn implies that 
$\|P\|_{\logl\to\logl}\le \frac1\pi\Big\|\frac{1}{\sqrt{1-x^2}}\Big\|_{\logl}$.
Since $(1/\pi)\int_{-1}^11/\sqrt{1-x^2}\,dx=1$, it follows that $P$  is  
a linear projection from $\logl $ onto the one-dimensional subspace 
of $\logl$ spanned by $1/\sqrt{1-x^2}$.

Let $f\in\logl$. It follows from Theorem \ref{t-3-4} that 
$$
(TP)(f)=T(P(f))=\Big(\frac1\pi\int_{-1}^{1}f(t)\,dt\Big)T\Big(\frac{1}{\sqrt{1-t^2}}\Big)=0.
$$
Accordingly, $TP=0$.

(iv) Let $ f \in \logl$.  Since $\sqrt{1-t^2}\in L^\infty$,  the 
Poincare-Bertrand formula (cf.\ Corollary \ref{c-3-3}) yields 
$$
T\Big(\sqrt{1-t^2}\cdot T(f)(t) 
+f(t) T\big(\sqrt{1-s^2}\big)(t) \Big )(x) = T\big(\sqrt{1-t^2}\big)(x)T(f)(x) - \sqrt{1-x^2}f(x)
$$ 
for  a.e. $x \in  (-1,1)$ . So,  recalling that $T\big(\sqrt{1-s^2}\big) (t) 
= -t$ for $t \in (-1,1)$  (via  \cite[\S 4.3 (9)]{tricomi} or 
\cite[(11.57)]{king}) the previous identity gives
\begin{equation}\label{eq-4-7}
T\big(\sqrt{1-t^2} T(f) -tf(t)  \big) (x)= -xT(f)(x)-\sqrt{1-x^2}f(x),
\quad\text{a.e. } |x|<1,
\end{equation}
that is, as an equality of functions in $L^0$.  From the proof of Theorem \ref{t-3-4}
we have that
$$
T(tf(t))(x)= \frac{1}{\pi}\int_{-1}^1f(t)\, dt + x T(f)(x),
\quad\text{a.e. $x \in (-1,1)$.}
$$
This and (\ref{eq-4-7}) give
$$
T\big(\sqrt{1-t^2} T(f) (t) \big)(x)  -  \frac{1}{\pi}\int_{-1}^1f(t)\, dt = -\sqrt{1-x^2}f(x)
$$
for a.e. $x \in (-1,1)$ and hence,
\begin{align*}
-\widehat{T}(T(f))(x)&=
\frac{1}{\sqrt{1-x^2}}T\big(\sqrt{1-t^2} T(f) (t) \big)(x)  
\\ &=
\frac{1}{ \sqrt{1-x^2}} \left(\frac{1}{\pi}\int_{-1}^1f(t)\, dt\right)  -f(x) 
= P(f)(x) - f(x)
\end{align*}
for a.e. $x \in (-1,1)$. This establishes \eqref{eq-4-5} as  an identity in $L^0$.

Since $T(f)\in L^1$ and $\widehat{T}\colon L^1\to L^0$, we have that
$\widehat{T}(T(f))\in L^0$. Because $f-P(f)\in\logl$, it follows from 
\eqref{eq-4-5} that $\widehat{T}(T(f))\in \logl$.
So, $\widehat{T}T\colon \logl\to L^0$ actually takes its values in $\logl$
and, being equal to $I-P$ (cf.  \eqref{eq-4-5}), it is continuous from $\logl$ into itself.
The estimate \eqref{eq-4-6} follows from (iii) and \eqref{eq-4-5}.
\end{proof}


An important consequence of Theorem \ref{t-4-10} is the following
description of the range space $T(\logl)$.

\begin{corollary}\label{cor-4-10}  
A function $g\in L^1$ belongs to the range space
$T(\logl)$ if and only if it satisfies  both $\widehat T(g)\in \logl$ and
$T(\widehat{T}(g))=g$. That is, 
$$
T(\logl)=\Big\{g\in L^1:  \widehat T(g)\in \logl,\; T(\widehat{T}(g))=g\Big\}.
$$
\end{corollary}

\begin{proof}
Let $g\in L^1$ satisfy $\widehat T(g)\in \logl$ and $T(\widehat T(g))=g$. Then
$g\in T(\logl)$.

Conversely, let $g\in T(\logl)$. Then there exists $h\in \logl$ such that $g=T(h)$. 
Theorem \ref{t-4-10}(iv) yields
$$
\widehat{T}(g)=\widehat{T}(T(h))=(I-P)(h)\in \logl.
$$
From Theorem \ref{t-4-10}(iii)  it follows that $T(P(h))=0$  and so
$$
T(\widehat{T}(g))=T(h-P(h))=T(h)-T(P(h))=T(h)=g.
$$
\end{proof}


The description of the range space $T(\logl)$ obtained 
in  Corollary \ref{cor-4-10}, although precise and useful 
(as we will soon see),
is not a full identification. The next result presents various facts
aimed at a better understanding  of  $T(\logl)$.  We highlight statement (ii)
which shows that $\logld\subseteq T(\logl)$.

\begin{proposition}\label{p-4-11} 
The following assertions hold for the continuous linear operators  $\tlog\colon \logl \to L^1$
and $\widehat T\colon \logl \to L^1$.
\begin{itemize}
\item[(i)] The range  $T(\logl)$ is a proper dense linear subspace of $L^1$.
\item[(ii)]  $\logld$ is included in the range $T(\logl)$.
\item[(iii)] The range $T(\logl)$  is not included in $\logl$.
\item[(iv)] The range  $\widehat T(\logl)$  is not included in $\logl$.
\item[(v)]   $\logl$ is not included in the range $T(\logl)$.
\item[(vi)]  $\logl$ is not included in the range $\widehat T(\logl)$.
\end{itemize}
\end{proposition}


\begin{proof}
(i)  In order to  exhibit a function $g\in L^1$ which is not in  $T(\logl)$,  
define  the  integrable  function $h$ on $(-1,1) $ by $h(t) := 1/\big({t(\log t)^2}\big)$ 
for $t \in (0,1/2)$ and $h(t):=0$ for $t \in (-1,1)\setminus(0,1/2)$.  It follows from 
Theorem 1(b) and its proof in \cite{kober}  that $ T(h)$ is not integrable over the 
subinterval $(-1/2,0)$.   So, also  $T(h)(x)/\sqrt{1-x^2}$   is not integrable 
over $(-1/2,0)$ since
$1 < 1/\sqrt{1-x^2} < 2/\sqrt{3}$ for $x \in (-1/2, 0)$.  
Hence, $T(h)(x)/\sqrt{1-x^2}\not\in L^1$.
Now, observe that  the function  $g(t) : = h(t)/\sqrt{1-t^2}$ belongs to $L^1$.  
To show that   $g\not\in T(\logl)$, assume, on the contrary, 
that $ g= T(f)$ for some $ f \in \logl$.  
Applying   Corollary \ref{cor-4-10}  yields $\widehat T(g)\in\logl$, where
$$
\widehat T(g)(x)=\frac{-1}{\sqrt{1-x^2}}  T\big(\sqrt{1-t^2}\,g(t)\big)(x)
= \frac{-T(h)(x)}{\sqrt{1-x^2}} , \quad x \in (-1,1);
$$
see \eqref{eq-4-1}. This implies that $T(h)(x)/\sqrt{1-x^2}\in \logl \subseteq L^1$;
contradiction.  So, 
$g \in L^1\setminus  T(\logl)$ and hence, 
$T(\logl)$ is strictly smaller than $L^1$.

Since $\bigcup_{1<p\le\infty}L^p\subseteq \logl$,   it follows
from Theorem \ref{t-4-10}(ii)   that 
$T(\widehat T(g))=g$ for every $g\in\bigcup_{1<p<2}L^p=\bigcup_{1<p\le\infty}L^p$.
Moreover, if $g\in L^p$ for some $1<p<2$,
then $\widehat T(g)\in L^p\subseteq\logl$. Accordingly,
Corollary \ref{cor-4-10} implies that
$\bigcup_{1<p\le\infty}L^p\subseteq T(\logl)$ and so 
$T(\logl)$ is dense in $L^1$.

(ii) Let $g\in\logld$. It follows from
Proposition \ref{p-4-5} that $\widehat T(g)\in\logl$. Moreover,
since $g\in\logl$, it follows from
Theorem \ref{t-4-10}(ii) that $\widehat T(T(g))=g$. 
Then    Corollary \ref{cor-4-10} implies that $g\in T(\logl)$.

(iii) Suppose that $T(\logl) \subseteq \logl$.   
Then $T\colon\logl\to\logl$.  To apply the closed graph theorem, consider a 
sequence $(f_n)_{n=1}^\infty$ such that  $f_n \to 0$ in $\logl$ and $T (f_n) \to g$ 
in $\logl$.  Theorem \ref{t-2-1} implies that  
$T(f_n ) \to 0$ in $L^1$.  On the other 
hand, the continuous inclusion $\logl \subseteq L^1$ ensures  that 
$T (f_n) \to g$ in $L^1$.  So, $g= 0$  in $L^1$ and 
hence, also in $\logl$.  Accordingly,
$T\colon \logl\to\logl$ is a closed operator and hence, it is continuous. This is impossible  
because both of the Boyd indices of $\logl$ are $1$.

(iv) Assume, by way of contradiction, that $\widehat T(\logl)\subseteq \logl$.
Let $f\in\logl$ be any function supported   in $(-1/2,1/2)$. The claim is that 
$$
\chi_{(-1/2,1/2)}T(f)\in \logl.
$$ 
To see this, note that
$\sqrt3/2<\sqrt{1-x^2}<1$ for $x\in (-1/2,1/2)$. Then $f(x)/\sqrt{1-x^2}$ is supported in $(-1/2,1/2)$ and
belongs to $\logl$. Our assumption implies, via \eqref{eq-4-1}, that 
$$
\widehat T\Big(\frac{f(t)}{\sqrt{1-t^2}}\Big)(x)=\frac{-1}{\sqrt{1-x^2}}\, T(f)(x)\in \logl.
$$
Consequently
$$
\chi_{(-1/2,1/2)}(x)T(f)(x)=\Big(\chi_{(-1/2,1/2)}(x)\sqrt{1-x^2}\Big)\frac{1}{\sqrt{1-x^2}}T(f)(x)\in \logl.
$$
What has been proved implies that the map
$$
T\colon \logl(-1/2,1/2)\to \logl(-1/2,1/2)
$$ 
is well defined. But this cannot be since $\logl(-1/2,1/2)$ has trivial Boyd indices. It follows that
 $\widehat T(\logl)\not\subseteq \logl$.

(v) Assume, by way of contradiction, that $\logl \subseteq T(\logl)$.
The claim is, in this case, that $\widehat T(\logl)\subseteq \logl$.
To see this let $g\in\logl$. By assumption,  $g\in T(\logl)$. So,
there exists $f\in\logl$ satisfying $T(f)=g$. It follows
from \eqref{eq-4-5} that  
$$
\widehat T(g)=f-P(f)\in \logl.
$$ 
Hence, $\widehat T(\logl)\subseteq \logl$ which contradicts part (iv).

(vi)  Assume, by way of contradiction, that $\logl \subseteq \widehat T(\logl)$. 
Then   $T(\logl ) \subseteq T\big(\widehat T(\logl)\big)$.  
Theorem \ref{t-4-10}(ii) implies that $T(\widehat T(\logl))= \logl$.
But, then $T(\logl ) \subseteq T\big(\widehat T(\logl)\big) = \logl$, 
which contradicts part (iii).
\end{proof}


\begin{remark}\label{r-4-12}
(i) Proposition \ref{p-4-11}(i) implies that
$T(\logl)$  is not closed in $L^1$. Consequently, the operator $\tlog$ is not Fredholm.

(ii) As noted in the proof of Proposition \ref{p-4-11}(i) we have  
that $\bigcup_{1<p\le\infty} L^p\subseteq  T(\logl)$.
  The inclusion is actually \textit{proper}. To see this, let 
$g\in \logl\setminus \bigcup_{1<p\le\infty} L^p$. If we had equality,
then $T(g)\in L^q$ for some $q\in(1,2)$ and hence, $\widehat T(T(g))\in L^q$.
Since also $1/\sqrt{1-x^2}\in L^q$, it follows that 
$\widehat T(T(g))+c/\sqrt{1-x^2}$ belongs to $L^q$ for all $c\in\C$.
Then \eqref{eq-4-5} implies that $g\in L^q$, which contradicts the choice of
$g$. So $\bigcup_{1<p\le\infty} L^p\subsetneqq  T(\logl)$.
\end{remark}


Recall that the airfoil equation is $T(f)=g$, that is, 
\begin{equation}\label{eq-4-9}
\frac{1}{\pi}\, \mathrm{ p.v.}  \int_{-1}^{1}\frac{f(x)}{x-t}\,dx=g(t),
\quad \mathrm{a.e. }\; t\in(-1,1).
\end{equation}
Given $g$, an  \textit{inversion formula} is needed to solve this equation for $f$. 
We are now able to  obtain such an inversion formula  for 
 $\tlog\colon\logl\to L^1$, via  the extended Poincar\`{e}-Bertrand formula in Corollary \ref{c-3-3}
and the properties of $\tlog$ presented in Theorem \ref{t-4-10}. 
and Corollary \ref{cor-4-10}.


\begin{theorem}\label{t-4-13}
Let $g$ belong to the range of  $\tlog\colon\logl\to L^1$.  
Then all solutions  $f\in \logl$ of the airfoil equation \eqref{eq-4-9} are 
of the form $\widehat{T}(g)(x)+c(1-x^2)^{-1/2}$, for $|x|<1$, that is,  
\begin{equation}\label{eq-4-10}
f(x)=\frac{-1}{\sqrt{1-x^2}}\; T\left(\sqrt{1-t^2}\, g(t)\right) (x)
+ \frac{c}{\sqrt{1-x^2}},\quad \mathrm{a.e. } \;|x|<1,
\end{equation}
with $c\in\mathbb{C}$ arbitrary.  In particular,
\eqref{eq-4-10} is satisfied for   every $g\in \logld$.
\end{theorem}

\begin{proof}
Given $c\in\mathbb{C}$ define $f$ by \eqref{eq-4-10}. Since $g\in T(\logl)$,
Corollary \ref{cor-4-10} implies that $\widehat T(g)\in\logl$.
Since $f=\widehat T(g)+c/\sqrt{1-x^2}$, it follows that $f\in\logl$.
From Theorem \ref{t-3-4} we have $T(c/\sqrt{1-x^2})=0$ and so 
$T(f)=T(\widehat{T}(g))$. Again
Corollary \ref{cor-4-10}   gives $T(\widehat{T}(g))=g$
  and so $f$ is a solution of the airfoil equation.

Conversely, let $f\in \logl$ satisfy   $T(f)=g$. 
It follows from Theorem \ref{t-4-10}(iv)   that  
$$
f=\widehat{T}(T(f))+P(f)=\widehat{T}(g)+P(f),
$$ 
with $P(f)=c/\sqrt{1-x^2}$, where $c= \frac1\pi\int_{-1}^{1}f(t)\,dt$. 
Accordingly, $f$ has the form  \eqref{eq-4-10}.

If $g\in\logld$, then $g\in T(\logl)$ by Proposition \ref{p-4-11}(ii)
and so \eqref{eq-4-10} holds.
\end{proof}


A  consequence for $\tlog$, via   properties  of the range space
$T(\logl)$, is the following one.


\begin{proposition}\label{p-4-14}
The operator $\tlog\colon \logl \to L^1$ is not  compact. 
\end{proposition}

\begin{proof} 
Consider the set $W:=\{\chi_A:A\in\mathcal{B}\}$. It is a bounded subset of $L^\infty$ 
and hence,  also a bounded subset of $\logld$.
Since $\widehat T\colon\logld\to\logl$ continuously (cf.\ Proposition \ref{p-4-5}),
it follows that $\widehat T(W)=\{\widehat T(\chi_A):A\in\mathcal{B}\}$ is a bounded 
subset of $\logl$.   Moreover,  $\chi_A\in\logl$   for each $A\in\mathcal{B}$ and so 
Theorem \ref{t-4-10}(ii)   implies  that $T(\widehat T(\chi_A))=\chi_A$. 
Thus, $W=T(\widehat T(W))\subseteq T(\logl)\subseteq L^1$.
But, $W$ is a well known example of a  subset of $L^1$
which is not relatively compact; see \cite[Example III.1.2]{diestel}. Hence,  
$\tlog\colon \logl \to L^1$ is not a compact operator.
\end{proof}


\section{The optimal domain for $T$ taking values in $L^1$}
\label{S5}


In this section we  prove that  the operator $\tlog\colon\logl \to L^1$  is optimally defined. 
The approach used for this   in the case when $T\colon X\to X$, with $X$
a r.i.\ space having non-trivial Boyd indices,  
was to explicitly determine
the  \textit{optimal domain} of $T$ (with values in $X$),
denoted by $[T,X]$ (cf.\ \cite[\S3]{curbera-okada-ricker-mh}), 
and then to show
that $[T,X]$ is isomorphic with  $X$ itself. This strategy 
was rather involved and required an in-depth study of 
$[T,X]$; cf.\ \cite[Theorem 5.3]{curbera-okada-ricker-ampa}, 
\cite[Theorem]{curbera-okada-ricker-mh}. 
We will profit from that work in this section.


Define the  space of functions
\begin{equation}\label{eq-5-1}
[T, L^1] : = \Big\{f \in L^1: T(h) \in L^1 \text{ for all } |h|  \le |f|\Big\}
\end{equation}
and the associated functional
\begin{equation}\label{eq-5-2}
\|f\|_{[T, L^1] }:= \sup\Big\{\|T(h)\|_{L^1}: |h|\le|f|
\Big\},\quad f \in [T, L^1] .
\end{equation}
We refer to $[T, L^1]$ as the  \textit{optimal domain} of $T$ with values in $L^1$.
Theorem \ref{t-2-1} implies that the inclusion $\logl \subseteq [T, L^1] $ is valid.
It is required to show that the 
opposite containment  $[T, L^1]\subseteq  \logl $ also  holds.


To investigate the structure and properties of $[T, L^1]$, we begin  
with the following result, which is modelled on \cite[Proposition 4.1 
and Lemma 4.3]{curbera-okada-ricker-ampa} for 
the case when $X$ is a r.i.\ space  with non-trivial Boyd indices.


\begin{proposition}\label{p-5-1} 
Let  $f \in L^1$. 
\begin{itemize}
\item[(i)]
The following conditions are equivalent.
\begin{itemize}
\item[(a)]  $f\in[T,L^1]$.
\item[(b)] $\displaystyle \sup_{|h|\le|f|}\|T(h)\|_{L^1}<\infty.$
\item[(c)] $T(f\chi_A)\in L^1$ for every $A\in\mathcal{B}$.
\item[(d)] $\displaystyle \sup_{A\in\mathcal{B}}\|T(f\chi_A)\|_{L^1}<\infty.$
\item[(e)] $T(\theta f)\in L^1$ for every $\theta\in L^\infty$ with $|\theta|=1$ a.e.
\item[(f)] $\displaystyle \sup_{|\theta|=1}\|T(\theta f)\|_{L^1}<\infty.$
\item[(g)] $fT(g) \in L^1$ for every $ g \in L^\infty$.
\end{itemize}
\item[(ii)]
Every function $f\in[T,L^1]$ satisfies the inequalities
\begin{equation*}
\sup_{A\in\mathcal{B}}\big\|T(\chi_A f)\big\|_X
\le
\sup_{|\theta|=1}\big\|T(\theta f)\big\|_{L^1}
\le
\sup_{|h|\le|f|}\big\|T(h)\big\|_{L^1}
\le
4 \sup_{A\in\mathcal{B}}\big\|T(\chi_A f)\big\|_{L^1}.
\end{equation*}
\item[(iii)] 
For every function $f\in[T,L^1]$ the following Parseval formula holds:
\begin{equation}\label{eq-5-3}
 \int_{-1}^1 fT(g) = - \int_{-1}^1 gT(f),\quad g\in L^\infty.
\end{equation}
\end{itemize}
\end{proposition}


\begin{proof}
(i) To establish that the conditions (a)--(f) are equivalent,
we can adapt the proof of \cite[Proposition 4.1]{curbera-okada-ricker-ampa}
to the case when $X=L^1$, after noting that the non-triviality of the Boyd indices
of $X$ is not used in that part of the proof.

(b)  $\Rightarrow$  (g).  Fix $g\in L^\infty$. 
Given $n\in\N$ define $A_n:=|f|^{-1}([0,n])$ and set
$f_n:=f\chi_{A_n}\in L^\infty$.
Since   $|f_n|\uparrow|f|$ pointwise on $(-1,1)$,
monotone convergence  yields 
\begin{equation}\label{eq-5-4}
\int_{-1}^1|f|\, |T(g)|=\lim_{n}\int_{-1}^1|f_n|\,|T(g)|.
\end{equation} 
Select $\theta_1,\theta_2\in L^\infty$ with $|\theta_1|=1$ and $|\theta_2|=1$ pointwise such 
that $|f|=\theta_1 f$ and $|T(g)|=\theta_2 T(g)$ pointwise. Then also
$|f_n|=\theta_1f_n$ pointwise for all $n\in\N$. Fix $n\in\N$. Then
$$
\int_{-1}^1|f_n|\,|T(g)|=\int_{-1}^1\theta_1\theta_2f_n\, T(g).
$$
Noting that both of the bounded functions $\theta_1\theta_2f_n$ and $g$ 
belong to $L^2$, we can apply the Parseval formula for $L^2$ to yield
$$
\int_{-1}^1|f_n|\,|T(g)|=-\int_{-1}^1gT(\theta_1\theta_2f_n).
$$
Since $\theta_1\theta_2f_n\in L^\infty$ implies that $T(\theta_1\theta_2f_n)\in L^1$,
it follows that 
$$
-\int_{-1}^1gT(\theta_1\theta_2f_n) \le \|g\|_\infty \|T(\theta_1\theta_2f_n) \|_{L^1}
\le \|g\|_\infty \sup_{|h|\le|f|}\|T(h)\|_{L^1}.
$$ 
Then \eqref{eq-5-4} implies that
$$
\int_{-1}^1|f|\,|T(g)|\le  \|g\|_\infty \sup_{|h|\le|f|}\|T(h)\|_{L^1}<\infty,
$$
that is, $fT(g)\in L^1$.

(g) $\Rightarrow$   (c). Fix $ A \in \mathcal{B}$.    Since $f\in L^1$,   we can choose
a sequence $(f_n)_{n=1}^\infty$ in $L^\infty$ such that  $|f_n|\uparrow|f|$ pointwise on $(-1,1)$. 
It follows that $f_n\chi_A\to f\chi_A$
in $L^1$ and hence, that $T(f_n\chi_A)\to T( f\chi_A)$ in measure (via Kolmogorov's 
theorem).  So, we can assume that $T(f_n\chi_A)\to T( f\chi_A)$ pointwise a.e. on 
$(-1,1)$, by passing to a subsequence, if necessary.   Next let $B\in \mathcal{B}$.    
Apply Parseval's formula  for $L^2$, which is permissible because  both  $f_n\chi_A$ and 
$ \chi_B$ belong to  $ L^\infty\subseteq L^2$, gives
$$
\int_BT(f_n\chi_A)  =  \int_{-1}^1\chi_BT(f_n\chi_A) = - \int_{-1}^1 (f_n\chi_A )T(\chi_B).
$$
Since $fT (\chi_B) \in L^1$ (by condition (g) with $g:= \chi_B$),  dominated convergence gives
$$
\lim_{n} \int_{-1}^1\chi_BT (f_n\chi_A) 
= - \lim_{n} \int_{-1}^1 (f_n\chi_A)T(\chi_B) = - \int_{-1}^1 (f\chi_A)T(\chi_B),
$$
that is,  the sequence $\left(\int_B T(f_n\chi_A)\right)_{n=1}^\infty$
is convergent in $\C$.  Noting that $B\in \mathcal{B}$ is arbitrary and 
$T(f_n\chi_A)\to T( f\chi_A)$ pointwise a.e.\ on $(-1,1)$,
this implies that $T(f_n\chi_A)\to T( f\chi_A)$ in $L^1$ and, in particular,
that $T(f\chi_A) \in L^1$.   So,  (c) is established.

Part (ii) can be proved as for  (4.1) in \cite{curbera-okada-ricker-ampa}, 
now with $X=L^1$; see p.1846 in \cite{curbera-okada-ricker-ampa}.

To prove part (iii), fix $g\in L^\infty$. Recall first that  $gT(f) \in L^1$  by (a)  
and that $fT(g) \in L^1$ by (g).  
Choose a sequence $(f_n)_{n=1}^\infty$ in $L^\infty$ such that  $|f_n|\uparrow|f|$ pointwise on $(-1,1)$. 
Parseval's formula for $L^2$ applied to $ f_n,\, g \in L^\infty\subseteq L^2$ 
gives $\int_{-1}^1 f_n T(g) = - \int_{-1}^1 g T(f_n)$, for every $n\in\N$.
Now, $\lim_n\int_{-1}^1 f_n T(g)= - \int_{-1}^1 f T(g)$ due
to dominated convergence as  $fT(g)\in L^1$.  Moreover, 
$T (f_n )\to T(f)$ in the norm of  $L^1$, which was 
established  in the proof of  (g) $\Rightarrow$ (c)  in part (i) 
(consider the case $A: = (-1,1)$ there).  In particular,
$\lim_n\int_{-1}^1 gT (f_n )=  \int_{-1}^1 g T(f)$.  
Combining these facts gives \eqref{eq-5-3} as
$$
\int_{-1}^1 fT(g) = \lim_{n}  \int_{-1}^1 f_n T(g) 
= -\lim_{n} \int_{-1}^1 g T(f_n) =- \int_{-1}^1 g T(f).
$$
\end{proof}


\begin{lemma}\label{l-5-2}
For every $f\in[T, L^1]$, both 
$f(x)\log (1-x)\in L^1$ and $f(x)\log (1+x)\in L^1$.
\end{lemma}

\begin{proof}
Let $f\in[T, L^1]$. According to Proposition \ref{p-5-1}(i) the function $fT(g)\in L^1$
for every $g\in L^\infty$. So, for $g=\chi_{(-1,1)}$, we have
$$
\frac1\pi f(x) \log\Big(\frac{1+x}{1-x}\Big)=(fT(g))(x)\in L^1.
$$
Moreover, the function
$$
\Big|\chi_{(0,1)}(x)f(x)\log(1-x)\Big|=
\Big|\chi_{(0,1)}(x)f(x)\log\Big(\frac{1+x}{1-x}\Big)-\chi_{(0,1)}(x)f(x)\log(1+x)\Big|
$$
belongs to $L^1$ because $\chi_{(0,1)}(x)f(x)\log(1+x)$ is bounded
on $(-1,1)$ as the singularity of $\log(1+x)$ is at $x=-1$. It follows that
$f(x)\log(1-x)\in  L^1$ because of the identity
$$
f(x)\log(1-x)=\chi_{(0,1)}(x)f(x)\log(1-x)+\chi_{(-1,0)}(x)f(x)\log(1-x),
\quad x\in(-1,1),
$$
and the fact that $\chi_{(-1,0)}(x)f(x)\log(1-x)$ is bounded 
on $(-1,1)$ as the singularity of $\log(1-x)$ is at $x=1$.

Similarly we can prove that $f(x)\log(1+x)\in  L^1$.
\end{proof}


We can now establish the following result alluded to above.

\begin{proposition}\label{p-5-3}
The space of functions $[T,L^1]\subseteq\logl$.
\end{proposition}

\begin{proof}
Let $f\in [T, L^1]$. Clearly \eqref{eq-5-1} implies that  $|f|\in[T, L^1]$ 
and so we can assume  that $f\ge0$. Of
course, also  $T(f)\in L^1$.
We use  a result of Stein for the space $\logl$, namely Theorem 3(b) 
in \cite{stein} for the case $n=1$ (i.e., in $\R$).
In this setting it is routine to check that the non-periodic Riesz
transform $R_1$ given in \cite[p.309]{stein} is precisely $-\pi T$
(i.e., the constant $c=\pi$). In order to deduce that
$f\in\logl$  from Theorem 3(b) in \cite{stein},
by selecting there $B_1=(-1,1)$ and $B_2=(-2,2)$,
we need to verify that
\begin{equation}\label{eq-5-5}
\int_{-2}^2|T(f)(x)|\,dx<\infty.
\end{equation}
By hypothesis $\int_{-1}^1|T(f)(x)|\,dx<\infty$  and so it remains to establish that
both $\int_{-2}^{-1}|T(f)(x)|\,dx<\infty$ and $\int_{1}^2|T(f)(x)|\,dx<\infty$. By an application
of Fubini's theorem, keeping in mind that $f\ge0$ and $(y-x)\ge0$ for 
$x\in(-2,-1)$ and $|y|<1$,  we have that
\begin{align*}
\int_{-2}^{-1}|T(f)(x)|\,dx&
=\int_{-2}^{-1} \frac1\pi \text{ p.v.} \int_{-1}^{1}\frac{f(y)}{y-x}\,dy\,dx
\\ &
=\frac1\pi \int_{-1}^{1}f(y)\left(\text{ p.v.} \int_{-2}^{-1}\frac{dx}{y-x}\right)\,dy
\\ &
=\frac1\pi \int_{-1}^{1}f(y)\log\left(\frac{y+2}{y+1}\right)\,dy
\\ &
=\frac1\pi \int_{-1}^{1}f(y)\log(y+2)\,dy-
\frac1\pi \int_{-1}^{1}f(y)\log(y+1)\,dy.
\end{align*}
Since $\log(y+2)$ is bounded on $(-1,1)$, as its singularity is at
$y=-2$, the first integral is finite. That the second integral is also finite
follows from Lemma \ref{l-5-2}. Consequently, $\int_{-2}^{-1}|T(f)(x)|\,dx<\infty$.

In a similar way it can be shown that  $\int_{1}^2|T(f)(x)|\,dx<\infty$. 
Accordingly, \eqref{eq-5-5} is satisfied and so $f\in\logl$. 
\end{proof}


The following two lemmas were established in 
Section 4 of \cite{curbera-okada-ricker-ampa} and
Section 2 of \cite{curbera-okada-ricker-qm}
for the optimal domain $[T,X]$ whenever  $X$ is a r.i.\ space having non-trivial
Boyd indices. The proofs given there are still valid in our setting, by adapting when necessary, 
the results used for $T\colon X\to X$ with the corresponding ones already established 
in this paper for $\tlog\colon \logl\to L^1$, namely  Theorem \ref{t-3-4}   for the kernel of
$\tlog$ and  Proposition \ref{p-5-1}(iii) for the
Parseval type formula valid for $[T, L^1]$.


\begin{lemma}\label{l-5-4}
The space $[T, L^1] $ is a linear lattice for the a.e.\ pointwise order
and, equipped with the  norm \eqref{eq-5-2}, it is a B.f.s. 
The norm \eqref{eq-5-2} is also given by
\begin{equation*}
\|f\|_{[T, L^1] }= \sup\Big\{\|fT(g)\|_{L^1}: g \in L^\infty, \,
\|g\|_{L^\infty} \le1\Big\},\quad f \in [T, L^1] .
\end{equation*}
\end{lemma}


\begin{lemma}\label{l-5-5}
The B.f.s.\ $[T, L^1] $ is the largest B.f.s.\ within $L^0$, containing $\logl$,  
to which   $\tlog\colon\logl \to L^1$ admits an $L^1$-valued, continuous linear extension.
\end{lemma}


We can now answer the question regarding the optimal extension of 
$\tlog\colon\logl \to L^1$. Of course, the operator $\tlog$ is extended
to $[T, L^1]$ in the obvious way.

\begin{theorem}\label{t-5-6}
The identity $[T, L^1] = \logl$ holds as an order 
and bicontinuous isomorphism between B.f.s.' and hence, 
$\tlog\colon\logl \to L^1$ does not admit a continuous linear  extension to 
any strictly  larger B.f.s. within $ L^0$.
\end{theorem}

\begin{proof}
Theorem \ref{t-2-1} and Proposition \ref{p-5-3} show that $\logl$ and
$[T, L^1]$ are equal as sets.

To prove that $\logl$ is continuously included in $[T, L^1] $,
let $f\in\logl$. For each $h\in L^0$ satisfying $|h|\le|f|$ we have that
$h\in\logl$ with $\|h\|_{\logl}\le\|f\|_{\logl}$. Theorem \ref{t-2-1} yields 
$\|T(h)\|_{L^1}\le \|T\|\cdot\|h\|_{\logl}\le \|T\|\cdot\|f\|_{\logl}$ and so we can
conclude from \eqref{eq-5-2} that 
$$
\|f\|_{[T, L^1]}\le \|T\|\cdot\|f\|_{\logl},\quad f\in\logl,
$$
where $\|T\|$ denotes $\|T\|_{\logl\to L^1}$.
That is, the natural inclusion $\logl\subseteq[T, L^1]$ is continuous.

Since both $\logl$ and $[T, L^1]$ are B.f.s.' (cf. Lemma \ref{l-5-4}) it follows
that the natural inclusion $\logl\subseteq[T, L^1]$ is an isomorphism.

Lemma \ref{l-5-5} implies that $\tlog\colon\logl \to L^1$ 
does not admit a continuous linear  extension to any strictly  larger B.f.s.\  within $ L^0$.
\end{proof}


\begin{remark}\label{r-5-7}  
The identification $[T, L^1]=\logl$ in Theorem \ref{t-5-6} 
implies that each of the equivalent conditions in Proposition \ref{p-5-1}(i)
\textit{characterizes} when a function $f\in L^1$ actually belongs to  $\logl$.  
In particular, for $f\in L^1$ we have that 
$$
f\in\logl \quad \text{if and only if}\quad  T(f\chi_A)\in L^1, \text{ for all } A\in\mathcal B.
$$
Since $\logl=[T,L^1]$, this is the same description arising for the optimal
domain $[F,\ell^{p'}(\Z)]$; see the question of R.~E.~Edwards
mentioned in the Introduction.
\end{remark}


It should be noted that there are functions in
$L^1\setminus \logl$ (i.e., outside of the domain of $\tlog$) which
$T$ nevertheless maps into $L^1$ (even into $\logl$).


\begin{proposition}\label{p-5-8}  
There exists $h\in L^1\setminus \logl$ such that
$T(h)\in\logl\setminus T(\logl)$.
\end{proposition}

\begin{proof} 
According to  Proposition \ref{p-4-11}(v) there exists $g\in\logl\setminus T(\logl)$, 
in which case $T(g)\in L^1$.
From Theorem \ref{t-4-10}(ii) it follows that $T(\widehat T(g))=g$.
The claim is that
\begin{equation}\label{eq-5-6}
\widehat T(g)\not\in\logl.
\end{equation}
Suppose not. Then $\widehat T(g)\in\logl$.   Together with $g=T(\widehat T(g))$   this 
implies that  $g\in T(\logl)$.  But, this contradicts the   choice of $g$. So,  \eqref{eq-5-6} holds.

Set $h:=\widehat T(g)$. Then, $h\in L^1$ (as $g\in\logl$) and, from \eqref{eq-5-6},
we have that  $h\not\in \logl$. Moreover, 
$$
T(h)=T(\widehat T(g))=g\in\logl\setminus T(\logl).
$$
\end{proof}


\begin{remark}\label{r-5-9}
Proposition \ref{p-5-8} implies that the   linear subspace
$\{f\in L^1:T(f)\in L^1\}$ of $L^1$ is strictly larger than $\logl$.
However,   by Theorem \ref{t-5-6}  there must exist $g\in L^1$ with
$|g|\le|h|$, whenever   $h$ is a function as given  in Proposition \ref{p-5-8}, such that
$T(g)\notin L^1$.  The difference between the linear space 
$\{f\in L^1:T(f)\in L^1\}$ and the optimal domain $[T, L^1]$ ($=\logl$)
is that the latter  (as $\logl$ itself) is a function lattice. That is, it satisfies the ideal property 
namely, $f\in [T, L^1]$ and $|g|\le|f|$ a.e.\ imply  
that $f\in [T, L^1]$, whereas the former is not.
\end{remark}


\section{Appendix}
\label{S6}


\noindent
\textit{Proof of Lemma \ref{l-4-6}.}
For the definition of $c(\beta,\gamma)$ we refer to \eqref{eq-4-2}.
Concerning \eqref{eq-4-3}, observe that
$$
\int_{-1}^\infty\frac{1}{|\xi|^\beta(\xi +1)^\gamma}d\xi 
= \int_{-1}^0\frac{1}{|\xi|^\beta(\xi +1)^\gamma}d\xi 
+\int_{0}^1\frac{1}{\xi^\beta(\xi +1)^\gamma}d\xi 
+ \int_1^\infty\frac{1}{\xi^\beta(\xi +1)^\gamma}d\xi.
$$

For the first integral note that
\begin{align*}
\int_{-1}^{0}\frac{d\xi}{|\xi|^\beta|\xi+1|^\gamma}
&=
\int_{0}^{1}\frac{d\xi}{\xi^\beta(1-\xi)^\gamma}
=\left(\int_{0}^{\frac12}+\int_{\frac12}^{1}\right)\frac{d\xi}{\xi^\beta(1-\xi)^\gamma}
\\ &\le
\int_{0}^{\frac12}\frac{2^\gamma d\xi}{\xi^\beta}
+\int_{0}^{\frac12}\frac{2^\beta d\xi}{\xi^\gamma}
\le
2^\gamma\int_{0}^{1}\frac{ d\xi}{\xi^\beta}
+2^\beta\int_{0}^{1}\frac{ d\xi}{\xi^\gamma}
\\& \le
\frac{2}{1-\beta}+\frac{2}{1-\gamma}
\le 4\, c(\beta,\gamma).
\end{align*}

Next we have that
$$
\int_{0}^1\frac{1}{\xi^\beta(\xi +1)^\gamma}d\xi  
\le \int_{0}^1\frac{1}{\xi^\beta}d\xi = \frac{1}{1-\beta}\le  c(\beta,\gamma).
$$

Finally
$$\int_1^\infty\frac{1}{\xi^\beta(\xi +1)^\gamma}d\xi \le 
\int_1^\infty\frac{1}{\xi^\beta\xi ^\gamma}d\xi = \int_1^\infty\frac{1}{\xi^{\beta +\gamma}}d\xi 
= \frac{1}{\beta +\gamma -1}\le  c(\beta,\gamma).
$$

Thus
$$
\int_{-1}^\infty\frac{1}{|\xi|^\beta(\xi +1)^\gamma}d\xi \le 6 c(\beta,\gamma).
$$

Concerning \eqref{eq-4-4} we have, for fixed $x\in(-1,1)$, that
\begin{align*}
&\int_{-1}^1\frac{1}{|t-x| ^\beta (1-t^2) ^\gamma}dt
 =\int_{-1}^1\frac{1}{|t-x| ^\beta (1+t)^\gamma(1-t)^\gamma}dt\\
&= \int_{-1}^0\frac{1}{|t-x| ^\beta (1+t)^\gamma(1-t)^\gamma}dt 
+  \int_{0}^1\frac{1}{|t-x| ^\beta (1+t)^\gamma(1-t)^\gamma}dt\\&\le
\int_{-1}^0\frac{1}{|t-x| ^\beta (1+t)^\gamma}dt +  \int_{0}^{1}
\frac{1}{|t-x| ^\beta (1-t)^\gamma}dt.
\end{align*}
With the substitution $t-x  = \xi(x+1)$, we have
\begin{align*}
& \int_{-1}^0\frac{1}{|t-x| ^\beta (1+t)^\gamma}dt =  
\int_{-1}^{-x/(1+x)}\frac{1+x}{|\xi|^\beta |1+x|^\beta (\xi +1)^\gamma (1+x)^\gamma}d\xi 
\\ &
= \frac{1}{(1+x)^{\beta + \gamma -1}} \int_{-1}^{-x/(1+x)}\frac{1}{|\xi|^\beta (\xi +1)^\gamma}d\xi
\\ &
\le \frac{1}{(1+x)^{\beta + \gamma -1}} \int_{-1}^{\infty}\frac{1}{|\xi|^\beta (\xi +1)^\gamma}d\xi.
\end{align*}

Similarly,  the substitution $t-x  = \xi(x-1)$ yields
\begin{align*}
& \int_{0}^1\frac{1}{|t-x| ^\beta (1-t)^\gamma}dt =  
\int_{x/(1-x)}^{-1}\frac{x-1}{|\xi|^\beta |x-1|^\beta (\xi +1)^\gamma (1-x)^\gamma}d\xi 
\\ &
\le \frac{1}{(1-x)^{\beta + \gamma -1}} \int_{-1}^{x/(1-x)}\frac{1}{|\xi|^\beta (\xi +1)^\gamma}d\xi
\\ &
\le \frac{1}{(1-x)^{\beta + \gamma -1}} \int_{-1}^{\infty}\frac{1}{|\xi|^\beta (\xi +1)^\gamma}d\xi.
\end{align*}
In view of \eqref{eq-4-3}, we have that 
\begin{align*}
\int_{-1}^1\frac{1}{|t-x| ^\beta (1-t^2) ^\gamma}dt
& \le \left( \frac{1}{(1+x)^{\beta + \gamma -1}}+
\frac{1}{(1-x)^{\beta + \gamma -1}}\right) 
\int_{-1}^{\infty}\frac{1}{|\xi|^\beta (\xi +1)^\gamma}d\xi
\\ & = 
\frac{(1-x)^{\beta + \gamma -1} + (1+x)^{\beta + 
\gamma -1}}{(1-x^2)^{\beta + \gamma -1}}
\int_{-1}^{\infty}\frac{1}{|\xi|^\beta (\xi +1)^\gamma}d\xi
\\ &
\le \frac{2^{\beta +\gamma}\cdot 6\,c(\beta,\gamma)}{(1-x^2)^{\beta + \gamma -1}}
\le \frac{24\, c(\beta,\gamma)}{(1-x^2)^{\beta + \gamma -1}}.
\end{align*}
\qed


\begin{lemma}\label{l-6-1}
Let $1<p<3/2$ and define $\delta:=\frac23\cdot \frac{(p-1)}{p}$.   Then
\begin{equation}\label{eq-6-1}
c\Big(\frac{1}{2}, \delta p + \frac{1}{2}\Big) <\frac{3}{p-1}
\end{equation}
and
\begin{equation}\label{eq-6-2}
c\Big (\frac{1}{2} ,\delta p' \Big) \le 6.
\end{equation}
\end{lemma}

\begin{proof}
Note that $\min\{\frac12, p-1\}=p-1$ and that 
$\frac12(p-1)<\delta p<(p-1)$, as $\delta p=\frac23(p-1)$.
Let us first prove (\ref{eq-6-1}).  Observe   that 
\begin{align*}
c\Big(\frac{1}{2}, \delta p + \frac{1}{2}\Big)& =\max\left\{ \frac{1}{1-\frac{1}{2}},
\frac{1}{1-(\delta p +\frac{1}{2})}, \frac{1}{\frac{1}{2} +(\delta p +\frac{1}{2}) -1}\right\} \\
& =\max\left\{2,
\frac{1}{\frac{1}{2}-\delta p}, \frac{1}{\delta p }\right\}
< \max\left\{2,
6, \frac{2}{p-1 }\right\}.
\end{align*}
Here, the last inequality holds  because
$\frac{1}{\frac{1}{2} -\delta p} <6$,
which is a consequence of  the calculation
$ \frac{1}{2} -\delta p= \frac{1}{2} -\frac{2}{3}(p-1) =\frac{7-4p}{6}> \frac{1}{6}$  
(due to the assumption that $1< p < 3/2$,
which implies that $(7-4p)>1$),
and also because   $\delta p >\frac{1}{2}(p-1)$.  Accordingly, 
$c\big(\frac{1}{2}, \delta p + \frac{1}{2}\big) < 3/(p-1)$ 
because $3/(p-1) > 6$ via $ 1 < p < 3/2$.  This yields (\ref{eq-6-1}) as
$$
\left(c\Big(\frac{1}{2}, \delta p + 
\frac{1}{2}\Big)\right)^{1/p}< \left(\frac{3}{p-1}\right)^{1/p}< \frac{3}{p-1}.
$$

Next, since 
$$
\delta p' =\frac{2(p-1)}{3p}
\cdot\frac{p}{p-1}=\frac{2}{3},
$$ 
the inequality (\ref{eq-6-2}) follows as
$$
c\Big (\frac{1}{2} ,\delta p' \Big)= \max\left\{ \frac{1}{1-\frac{1}{2}},
\frac{1}{1-\delta p'}, \frac{1}{\frac{1}{2} +\delta p' -1}\right\}=\max\{2,3, 6\}=6.
$$
\end{proof}



\end{document}